\documentclass[12pt,a4paper]{amsart}
\usepackage{amsmath}
\usepackage{latexsym}
\newtheorem{thm}{Theorem}
\newtheorem{prop}{Proposition}
\newtheorem{lemma}{Lemma}

\newtheorem{defi}{Definition}
\newtheorem{cor}{Corollary}
\newcommand{\mcl}{\mathcal}
\newcommand{\eps}{\varepsilon}
\newcommand{\Om}{\Omega}
\newcommand{\Ga}{\Gamma}

\newcommand{\Sgm}{\Sigma}
\newcommand{\KOmm}{$K\Omega$-{\it mod}}
\newcommand{\KGam}{$K\Ga$-{\it mod}}

\newcommand{\mKOm}{{\it mod\/}-$K\Omega$}

\newcommand{\mKGa}{{\it mod\/}-$K\Gamma$}
\newcommand{\Cdnga}{C^{\Gamma}}
\newcommand{\Rdnga}{R^{\Gamma}}

\newcommand{\Cupga}{C_{\Gamma}}
\newcommand{\Rupga}{R_{\Gamma}}

\newcommand{\Hom}{\operatorname{Hom}}
\newcommand{\Iso}{\operatorname{Iso}}
\newcommand{\St}{\operatorname{St}}
\newcommand{\sk}{\operatorname{sk}}
\newcommand{\stk}{\stackrel}
\newcommand{\til}{\tilde}
\newcommand{\ulin}{\underline}
\newcommand{\olin}{\overline}

\newcommand{\zed}{\mathbf{Z}}
\newcommand{\lra}{\longrightarrow}
\newcommand{\lla}{\longleftarrow}
\newcommand{\Ext}{\operatorname{Ext}}
\newcommand{\Lie}{\operatorname{Lie}}
\newcommand{\Tor}{\operatorname{Tor}}
\newcommand{\Map}{\operatorname{Map}}
\newcommand{\Ind}{\operatorname{Ind}}
\newcommand{\Res}{\operatorname{Res}}
\newcommand{\from}{\leftarrow}
\newcommand{\coker}{\operatorname{coker}}
\newcommand{\cro}{\operatorname{cr}}
\newcommand{\crsh}{\cro^\sharp}
\newcommand{\del}{\partial}
\newcommand{\im}{\operatorname{im}}
\newcommand{\Ho}{\operatorname{Ho}}

\newcommand{\Calend}{\mcl{E}\!\:\!\mathit{nd}}

\begin{document}

\title{$E_\infty$ obstruction theory}
\author{Alan Robinson}
\address{Mathematics Institute, University of Warwick, Coventry CV4 7AL, United Kingdom}
\email{C.A.Robinson(at)warwick.ac.uk}
\begin{abstract} The space of $E_\infty$-structures on an simplicial operad~$\mcl{C}$ is the limit of a tower of fibrations, so its homotopy is the abutment of a Bousfield-Kan fringed spectral sequence.  The spectral sequence begins (under mild restrictions) with the stable cohomotopy of the right $\Ga$-module $\pi_*\mcl{C}$; the fringe contains an obstruction theory for the existence of $E_\infty$-structures on~$\mcl{C}$. This formulation is very flexible: applications extend beyond structures on classical ring spectra to examples in motivic homotopy theory. \end{abstract}
\date{January 6, 2013}
\maketitle
\vspace{-7mm}

\tableofcontents

\section{Introduction}

This paper develops obstruction theory for $E_\infty$ operads in more detail than the author's brief account in \S5 of  \cite{rob:inv}.  We study the classification of maps into a general simplicial (or topological) operad $\mcl{C}$ from a standard $E_\infty$ operad $\mcl{T}$ (which is in fact the operad of trees introduced in \cite{rob-whi2}).  One might see this as related to lifting results in~\cite{boa} and in Chapter III of Boardman and Vogt's classic work \cite{boa-vog}.  Our proofs are completely different from those in~\cite{rob:inv}.  The advantage of our approach here is that theorems apply to operads acting not only in category of spectra, but also in a much broader context.  The usefulness of such an approach became apparent recently in an application to motivic homotopy theory~\cite{nau-spi-ost}.

Our main results are proved in \S\ref{sec: e1} -- \S\ref{sec: fringe}.   Theorem~\ref{thm: e2} gives a Bousfield-Kan type fringed spectral sequence for the homotopy groups of the space of $E_\infty$ structures on a simplicial (or topological) operad  $\mcl{C}$ with a given homotopy associative and commutative pairing.  The $E_2$ term of the spectral sequence consists of certain $\Ext$ groups which by work of Pirashvili~\cite{pir:hodge} may be interpreted as stable cohomotopy groups of the $\Ga$-module of homotopy groups $\pi_*\mcl{C}$, which is described 
in~\ref{sec: piop}. 
Theorem~\ref{thm: obst} establishes the existence of a sequence of  obstructions to the existence of such an $E_\infty$ structure.  The obstructions lie in cohomotopy groups $\pi^{n}\pi_{n-2}\,\mcl{C}$ beyond the fringe of the aforementioned spectral sequence.  Similar situations have been investigated in~\cite{bousfield:1989}, in~\cite{goe-hop}, and in unpublished work by the same authors.

\section{Operads}

We work principally with operads of spaces.  By the term \emph{space} we here always mean \emph{simplicial set}, unless the context requires otherwise.

Our operads may not have units, so we adapt the definition of Markl \cite{Markl:1996} to the category of spaces.   As in \cite{rob-whi2} we index the spaces in an operad $\mcl{C}$ not by non-negative integers, but by all finite sets $S$ (in a suitable universe).  This eliminates from our proofs a host of spurious complexity which tends to arise because  our main constructions proceed by induction over the \emph{quotient} sets of $S$, which have no sufficiently standard enumeration.

If $\mcl{C}$ is an operad and $S$ happens to be the set $\{1,2,\dots,n\}$ then we write $\mcl{C}_n$ for $\mcl{C}_{\{1,2,\dots,n\}}$.   The operad $\mcl{C}$ is determined up to isomorphism by its values $\mcl{C}_n$ on this spine and the compositions among them.  The space $\mcl{C}_S$ has various right compositions $\circ_s$ with $\mcl{C}_T$ corresponding to differently labelled inputs $s\in S$.  The deleted disjoint union $(S\setminus\{s\})\sqcup T$ is denoted by $S\sqcup_s T$.   The formal definition is as follows.

\begin{defi} An \label{def: operad} \emph{operad} $\mcl{C}$ prescribes \begin{enumerate}
  \item For each finite set $S$ a space $\mcl{C}_S$, cofunctorial with respect to isomorphisms $\phi:S_1 \to S_2$ of finite sets
  \item  For each pair $(S,T)$ of finite sets and each $s\in S$ a natural composition map 
$\circ_s:\mcl{C}_S \times \mcl{C}_T \lra \mcl{C}_{S\sqcup_s T}$, such that two conditions hold:
  \item \emph{(Left-right associativity)} if $s\in S$ and $t\in T$ we have 
$$(x \circ_s y) \circ_t z \quad = \quad x \circ_s (y \circ_t z)$$   for all $x\in \mcl{C}_S, \: y\in \mcl{C}_T, \: z\in \mcl{C}_U$.
  \item \emph{(Right-right associativity)} if $s, s' \in S$ we have 
$$(x \circ_s y) \circ_{s'} z \quad = \quad (x \circ_{s'} z) \circ_s y$$ for $x,y,z$ as above.  \end{enumerate} \end{defi}

We now mention unit elements for operads.  The space $\mcl{C}_{\{s\}}$ corresponding to a singleton is determined up to unique isomorphism by clause (1) of the above definition, so we may identify $\mcl{C}_1$ with $\mcl{C}_{\{s\}}$ for any $s$.  A \emph{unit} for $\mcl{C}$ is an element $\eps$ of $\mcl{C}_1$ (or equivalently of $\mcl{C}_{\{s\}}$, where $s\in S$) such that for each finite set $S$ 
\begin{itemize} \item left composition with $\eps\in \mcl{C}_1$ is the identity map on $\mcl{C}_S$
             \item any right composition $\circ_s$ with $\eps\in \mcl{C}_{\{s\}}$ is the identity map on $\mcl{C}_S$
\end{itemize}
Strict units, when they exist, are unique.  It often suffices to have only homotopy units, which are defined by replacing ``the identity map'' by ``homotopic to the identity map'' in the above.

In $A_\infty$ theory there is a standard operad $\mcl{A}$ consisting of the Stasheff polytopes.  An action of $\mcl{A}$ is an $A_\infty$ structure.  We shall show that in $E_\infty$ theory there is likewise a canonical example with a central r\^ole. 
This is the operad $\mcl{T}$  of \cite{rob-whi2}, defined as the product $\mcl{P} \times \mcl{Q}$, where~$\mcl{P}$ is the partitions operad and $\mcl{Q}$ is the Barratt-Eccles operad \cite{bar-ecc1}.

\subsection{The operad $\mathcal{P}$ of partitions}

We define the \emph{partition operad} $\mcl{P}$ in terms of trees, and explain the connection with partitions later in the section.  So we assign to a finite set $S$ the space $\mcl{P}_S$ of trees \cite{rob-whi1} with internal edges of lengths between 0 and 1, and with twigs labelled by all elements of the set $S$.  As we are using operads without unit, we shall define $\mcl{P}_S$ to be empty in the trivial cases $|S| \le 1$.  (The alternative would be to admit as unit the unique tree with one twig.)  We summarise the properties which we shall need, and refer to \cite{rob-whi1} and \cite{rob:par} for proofs.

 The composition map
$$\circ_+\colon \mcl{P}_M \times \mcl{P}_{N^{\!+}}  \lra \mcl{P}_S \quad \text{where} \quad |M|,|N^{\!+}| \ge 2$$ grafts the root of $\mcl{P}_M$ into the twig labelled $+$ in $\mcl{P}_{N^{\!+}}$, forming a new internal edge of length 1.  
Thus a tree is a composite iff it has an internal edge of maximal length.  Such trees are also called \emph{decomposable} or \emph{fully grown}.   Each space $\mcl{P}_S$ with $|S| \ge 2$ can be contracted, by shrinking edge lengths of trees, to the point consisting of the \emph{star tree} with twigs $S$.

\begin{prop} \emph{(Properties of the operad of partitions.)} \label{prop: partop}
 \begin{enumerate}                                                        
                     \item Each composition map 
$$\mcl{P}_M \times \mcl{P}_{N^{\!+}}  \lra \mcl{P}_S$$ 
is injective: its image is the {\bf face} of $\mcl{P}_S$ corresponding to the partition $S = M \sqcup N$ of $S$.
                      \item Two different faces of $\mcl{P}_S$ are disjoint, or intersect in a common subface.  Their union is the \emph{boundary} $\partial\mcl{P}_S$, the subspace of all decomposable trees.
                      \item The space  $\mcl{P}_S$ is a cone, with the star tree as vertex and the boundary $\partial\mcl{P}_S$ as its base. $\Box$
\end{enumerate} \end{prop}

This type of face structure is the vital characteristic of \emph{cofibrant operads} as defined in \cite{rob-whi2}.  The term ``cofibrant'' is used informally: we shall not need any model structure on the category of operads.

The group of permutations of $S$ acts on $\mcl{P}_S$ on the right, by permuting the labels on the twigs.
Yet the operad $\mcl{P}$ {\bf fails} to be an $E_\infty$ operad, because this action is not free.  We shall remedy this by constructing the operad $\mcl{T}$ in \ref{sec: treeop}.

\begin{prop} \emph{(Homotopy type of $\partial\mcl{P}_S$.)} \label{prop: partbdry} Let $n = |S|$.
\begin{enumerate}
\item The boundary $\partial\mcl{P}_S$ (the union of all its faces) has the homotopy type of a wedge of spheres: in fact
$$\partial\mcl{P}_S \: \simeq \: \bigvee_{(n-1)!}S^{n-3} \qquad \text{and} \qquad 
\mcl{P}_S/{\partial\mcl{P}_S} \: \simeq \: \bigvee_{(n-1)!}S^{n-2}\; .$$

\item Under the action of the symmetric group $\Sigma_n$,  the cohomology 
$\til H^{n-3}(\partial\mcl{P}_S)$ is a $\zed\Sigma_n$-module isomorphic to $\epsilon\Lie_n$, the Lie representation twisted by the sign character.  $\Box$
\end{enumerate} \end{prop}

\begin{defi} \label{def: super} The representation $\eps\Lie_n$ which occurs in Proposition~\ref{prop: partbdry}  is called the \emph{$n$th superlie representation} $\mcl{S}_n$.  Its restriction to the subgroup $\Sigma_{n-1}$ is the regular representation.
Its $\zed$-dual, which is isomorphic to the homology group $\til H_{n-3}(\partial\mcl{P}_S)$, is denoted~$\mcl{S}^*_n$.
\end{defi}

The reason for the ``superlie'' terminology is explained in \cite{rob:par}; see also \cite{mel-reu, ree:1960}.
\bigskip

\subsection{Trees and partitions\label{sec: partree}}

We now reinterpret the spaces of the operad $\mcl{P}$, first in terms of partitions of $S$, and then in terms of a reduced nerve of a category of surjections under $S$. 

The set of partitions of the finite set $S$ is partially ordered.  We reverse the usual ordering by refinement, so that the discrete partition becomes the initial object, and the indiscrete partition the final object.  Let us assume $|S| \ge 2$. 
Partitions which are neither discrete nor indiscrete are called \emph{non-trivial}. 

\begin{lemma} \emph{(\cite{rob:par}, Prop. 2.7).} \label{lemma: treepart} The nerve of the partially ordered set of non-discrete partitions of $S$ is isomorphic to the space $\mcl{P}_S$ of $S$-trees. The nerve of the subset of non-trivial partitions corresponds to the space $\del\mcl{P}_S$ of decomposable trees.  $\Box$ \end{lemma}

The star-tree with twigs labelled by $S$ is a vertex corresponding to the trivial partition.  We have seen that $\mcl{P}_S$ is a cone with this star-tree as vertex.  The base of the cone is $\del\mcl{P}_S$.

For the next stage, note that partitions of $S$ correspond to equivalence relations, and therefore to quotient sets.  Taking a greater partition (in our ordering) corresponds precisely to quotienting further.

Let $\Om$ be the category of finite sets and surjective maps, and for $n \ge 0$ let $\ulin{n}$ be the object $\{1,2,\dots,n\}$ of $\Om$ .  We use the standard notation for over- and under-categories: thus when $S$ is a finite set, $S/\Om$ denotes the slice category of surjections of finite sets under~$S$.  Clearly 
$S/\Om \; \approx \; S/\Om/\ulin{1}\,$ since $\ulin{1}\,$ is a final object in $S/\Om$.

We also need some less standard notation for some full subcategories of the slice categories.  A surjection is \emph{strict} if it is not an isomorphism.

\begin{defi} \label{def: doubleslash} Let $S/\!/\Om$ be the subcategory of surjections of finite sets \emph{strictly} under $S$, that is, the full subcategory of $S/\Om$ whose objects are non-isomorphisms $S \to X$ in $\Om$. 

 Let $S{/\!/}\Om{/\!/}\underline{1}$ be the category of surjections of sets strictly under $S$ and strictly over the singleton object $\underline{1}$.  \end{defi}

\noindent Notice that the category $S{/\!/}\Om$ is obtained from $S{/\!/}\Om{/\!/}\ulin{1}$ by adding a final object, so that the nerve of $S{/\!/}\Om$ is the cone on the nerve of the subcategory.

\begin{lemma} \label{lemma: slicetree} There is a natural weak homotopy equivalence of pairs 
$$\Phi \colon \bigl(N\!(S{/\!/}\Om), N\!(S{/\!/}\Om{/\!/}\underline{1})\bigr)
 \quad \lra \quad \bigl(\mcl{P}_S, \del\mcl{P}_S\bigr)\; .$$
\end{lemma}

{\it Proof.} We map each object $S \to X$ of $S/\Om$  to the corresponding partition of $S$ by inverse images of points of $X$.   This maps the category $S{/\!/}\Om$ by an equivalence to the category of partitions of $S$ which are not discrete, and the subcategory $S{/\!/}\Om{/\!/}\underline{1}$ by an equivalence to the subcategory of partitions which are neither discrete nor indiscrete.  This gives a weak equivalence of pairs of nerves. Then we follow by the isomorphism of Lemma \ref{lemma: treepart} to the space of trees.  This gives the required homotopy equivalence. $\Box$

\bigskip

{\bf Notes.}  \begin{enumerate} \item One can restate the proof above by noting that in $S/\Om$ every endomorphism of an object is an identity map.  This implies that $\Phi$ is a map of nerves in which that inverse image of each point is contractible.  (Compare the last section of \cite{pir:dold}, which in turn refers to previous work of W.~L\"uck and of J.~Slominska.)

\item We can extend the above to describe the spaces in the operad $\mcl{P}$ in terms of nerves, but the description of the operad composition in terms of trees is more intuitive. 

\item The triangulation of $\mcl{P}_S$ arising from partitions, which we use here throughout, is finer than that used in \cite{rob-whi1}.

\item Proposition \ref{prop: partbdry} was proved in terms of partitions by M. L. Wachs after partial results by several workers \cite{wac:1998}.
\end{enumerate}
\medskip

\noindent In the following related proposition, $N$ denotes nerve, $C$ is cone and $\varSigma$~is suspension.

\begin{prop}\label{prop: nerpart} 
(i)  The nerve $N\!(S/\Om/\ulin{1})$ is weakly contractible.  \\
(ii)  Let $\del N\!(S/\Om/\ulin{1})$ denote the subcomplex of $N\!(S/\Om/\ulin{1})$ consisting of all simplices not having the 1-simplex $S\lra \ulin{1}$ as a face.  Then there is a natural weak equivalence of pairs
$$(N,\del N)(S/\Om/\ulin{1}) \quad \simeq \quad \bigl(C\varSigma(\del\mcl{P}_S),\varSigma(\del\mcl{P}_S)\bigr)\: .$$  
\end{prop}

{\it Proof.}  The category $S/\Om/\ulin{1}$ has an initial object $S$ and a terminal object $\ulin{1}$, so its nerve is weakly contractible.  Indeed the nerve can be regarded as a cone with either of these objects as base.  (Strictly speaking the initial object is $1_S$, \emph{etc.,} but we abuse notation here for simplicity.)

The simplices of $N(S/\Om/\ulin{1})$ not having $\ulin{1}\,$ as a vertex, form a subcone $C^+$ with vertex $S$ and base $N(S/\!/\Om/\!/\ulin{1})$.  The simplices of 
$N(S/\Om/\ulin{1})$ not having $S$ as a vertex, form a subcone $C^-$ with vertex $ \ulin{1}$ and the same base $N(S/\!/\Om/\!/\ulin{1})$.  The cones $C^+$ and $C^-$ intersect only in their common base, so their union $\del N\!(S/\Om/\ulin{1})$  is weakly equivalent to the suspension $\varSigma (N(S/\!/\Om/\!/\ulin{1}))$ and therefore (by Lemma~\ref{lemma: slicetree}) to $\varSigma(\del\mcl{P}_S)$.

Since by the first part $C^+ \cup C^-$ is (cofibrantly) embedded in the weakly contractible space $N\!(S/\Om/\ulin{1})$, we deduce the required weak equivalence of pairs.     $\Box$

\subsection{The Barratt-Eccles operad $\mcl{Q}$}

We now recall the \emph{Barratt-Eccles operad} $\mcl{Q}$.  In the original formulation \cite{bar-ecc1} its elements are the contractible free $\Sigma_n$-spaces $E\Sigma_n$.  In our categorical situation  operads have spaces indexed by finite sets $S$.  We define $\mcl{Q}_S$ to be the nerve of the slice category $S/\mcl{I}$ of \emph{isomorphisms} of finite sets under $S$.  Composition in the operad is deleted sum:  the map
$$\mcl{Q}_S \times \mcl{Q}_{T\!^+}  \lra  \mcl{Q}_{S\sqcup T}$$
takes the disjoint union of sets under $S$ with corresponding sets under $T$.  (The points which correspond under the isomorphisms to $+$ are deleted.)

The essential property of $\mcl{Q}$ is its freeness: different isomorphisms \mbox{$S' \to S$} in $\Om$ induce simplicial maps $\mcl{Q}_S \to \mcl{Q}_{S'}$ which are everywhere different.  (Equivalently, for every $n$ the symmetric group $\Sigma_n$ acts freely on $\mcl{Q}_{\underline{n}}$.)  Its disadvantage is that it is not a cofibrant operad: the decomposable elements do not behave as in Proposition \ref{prop: partop}.  Indeed a typical element of $\mcl{Q}$ is decomposable in many ways.

\subsection{The tree operad $\mcl{T}$\label{sec: treeop}}
 
This is our principal tool in investigating $E_\infty$ structures.  Again, it is an operad without unit in our formulation: the space $\mcl{T}_S$ is empty when $|S|\le 1$.  It is related to theories obtained by Boardman and Vogt's \hbox{$W$-construction} (\cite{boa-vog} Ch.III).

\begin{defi} The {\bf tree operad} $\mcl{T}$ is the product  $\mcl{P} \times \mcl{Q}$ of the Barratt-Eccles operad and the partition operad. \end{defi}

We can now investigate $\mcl{T}$ by using Propositions \ref{prop: partop} and \ref{prop: partbdry}.

\begin{prop} \emph{(Properties of the tree operad.)} \label{prop: treeop}
 \begin{enumerate}  \item  The tree operad  $\mcl{T}$ is $E_\infty$.                                                           
                     \item Each composition map 
$$\mcl{T}_M \times \mcl{T}_{N^{\!+}}  \lra \mcl{T}_S$$ 
is injective: its image is the {\bf face} of $\mcl{T}_S$ corresponding to the partition $S = M \sqcup N$ of $S$.
                      \item Two different faces of $\mcl{T}_S$ are disjoint, or intersect in a common subface.
                      \item The union $\partial\mcl{T}_S$ of all the faces (the {\bf boundary} of $\mcl{T}_S$) has the homotopy type of a wedge of spheres: in fact
$$\partial\mcl{T}_S \: \simeq \: \bigvee_{(n-1)!}S^{n-3} \qquad \text{and} \qquad 
\mcl{T}_S/{\partial\mcl{T}_S} \: \simeq \: \bigvee_{(n-1)!}S^{n-2}$$
where $n = |S|$.
                       \item Under the action of the symmetric group $\Sigma_n$,  the cohomology 
$\til H^{n-3}(\partial\mcl{T}_S)$ is a $\zed\Sigma_n$-module isomorphic to $\epsilon\Lie_n$, the Lie representation twisted by the sign character.  Its restriction to the subgroup $\Sigma_{n-1}$ is the regular representation. 
\end{enumerate} \end{prop}
                      
{\it Proof.  1.}  The space $\mcl{P}_S$ is contractible to the star tree with twigs $S$, by shrinking internal edges.  The category $S/\mcl{I}$ of isomorphisms of finite sets under $S$ has an initial object, so its nerve $\mcl{Q}_S$ is contractible.  Thus the product  $\mcl{T}_S = \mcl{P}_S \times \mcl{Q}_S$ is contractible for every $S$.  The freeness of the action of $\mcl{I}$ on $\mcl{T}$ follows from the freeness of the action on $\mcl{Q}$.  Therefore $\mcl{T}$ is $E_\infty$.  

2. The composition maps are injective both in the operad $\mcl{P}$ and in the operad $\mcl{Q}$.  Hence they are injective in the product $\mcl{T} = \mcl{P} \times \mcl{Q}$.

3. We first show that a (bi)simplex $(x,y) \in \mcl{P}_S \times \mcl{Q}_S = \mcl{T}_S$ lies on a face iff $x$ corresponds to a decomposable shape of $S$-tree, that is one with an internal edge of maximal length.  The necessity is obvious (from the definition of composition in $\mcl{P}$).  For sufficiency, suppose that the tree $x$ can be decomposed as $x_1 \circ x_2$, where $\circ$ denotes grafting of trees.  This gives a partition of $S$ as $M \sqcup N$, where $M$ (resp. $N\!^+$) labels the twigs of $x_1$ (resp. $x_2$).  Now we take the simplex 
$$y \quad = \quad \{S \lra S_0 \lra S_1 \lra \dots \lra S_k\}$$
of isomorphisms of sets under $S$, we partition each set $S_i = M_i \sqcup N_i$ and each isomorphism by transport from $S = M \sqcup N$, thus obtaining simplices 
$$y_1 \quad = \quad \{M \lra M_0 \lra M_1 \lra \dots \lra M_k\}$$ and (adjoining an extra input for grafting)
$$y_2 \quad = \quad \{N\!^+ \lra N_0\!^+ \lra N_1\!^+ \lra \dots \lra N_k\!^+\}$$
such that $(x.y)$ is the composite $(x_1,y_1) \circ (x_2,y_2)$, and therefore lies on a face.  We have incidentally shown that $\partial \mcl{T}_S \: = \: \partial\mcl{P}_S \times \mcl{Q}_S$.

Now suppose that $(x,y)$ lies on two distinct faces $F_1$ and $F_2$.  Then the tree $x$ has two maximal-length internal edges, along which it can be split into three parts.  Dissecting the label set $S$ into three much as above, we can decompose $y$ into three correspondingly.  We deduce that $(x,y)$ is a double composite.  Performing the compositions in two different orders, we deduce that $(x,y)$ lies in a face of $F_1$ which is also a face of $F_2$. 

4, 5.  The projection on the first factor
$$(\mcl{T}_S,\partial \mcl{T}_S) \: = \: (\mcl{P}_S,\partial\mcl{P}_S) \times \mcl{Q}_S \: \lra \: 
(\mcl{P}_S,\partial\mcl{P}_S)$$
is a homotopy equivalence of pairs, since $\mcl{Q}_S$ is contractible.  The results therefore follow immediately from Proposition \ref{prop: partbdry}.  This completes the proof of Proposition \ref{prop: treeop}.
$\Box$

\bigskip

Let $S = M \sqcup N^+$ be any non-trivial partition of $S$, as above.  We define the boundary of the face 
$\mcl{T}_M \times \mcl{T}_{N^+}$ of $\mcl{T}_S$ as usual to be 
$$\del(\mcl{T}_M \times \mcl{T}_{N^+}) \quad = \quad (\del\mcl{T}_M \times \mcl{T}_{N^+}) \cup 
(\mcl{T}_M \times \del\mcl{T}_{N^+})\;.$$
Thus $$(\mcl{T}_M \times \mcl{T}_{N^+})/\del{(\mcl{T}_M \times \mcl{T}_{N^+})} \quad \approx \quad 
(\mcl{T}_M/\del\mcl{T}_M) \wedge (\mcl{T}_{N^+}/\del\mcl{T}_{N^+})\;.$$
We denote by $\del\del\mcl{T}_S$ the union of the boundaries of all faces of~$\mcl{T}_S$.  The~next proposition follows from this and from statements {\it 2\/} and {\it 3\/} of Proposition \ref{prop: treeop}. 

\begin{prop} \label{prop: deldel} $$\begin{aligned} \del\mcl{T}_S/\del\del\mcl{T}_S \quad &\approx 
\quad \bigvee_{S = M\sqcup N} (\mcl{T}_M \times \mcl{T}_{N^+})/\del{(\mcl{T}_M \times \mcl{T}_{N^+})} \\
           &\approx \quad \bigvee_{S = M\sqcup N}\bigl((\mcl{T}_M/\del\mcl{T}_M) \wedge (\mcl{T}_{N^+}/\del\mcl{T}_{N^+})\bigr)\;. \end{aligned}$$
where the union is taken over all non-trivial partitions $S = M \sqcup N$ of $S$. $\Box$ \end{prop}

\bigskip

\section{Homological algebra of $\Gamma$-modules and $\Omega$-modules\label{sec: modules}}

\subsection{Indexing categories}

In the previous section we indexed the spaces in our operads by all finite sets (in some universe).  To handle the homology of our operads, we now develop some aspects of the homological algebra of modules over  categories of finite sets.  Much of the theory here is due to or developed from the ideas of Pirashvili \cite{pir:dold, pir:hodge} and Richter \cite{pir-richt}, and we refer to their original papers for some proofs.

\begin{defi}   We define two indexing categories as follows. \begin{enumerate}
\item[] $\Ga$ is the category of finite based sets and basepoint preserving maps.  
\item[] $\Om$ is the category of finite unbased sets and surjective maps. \end{enumerate} \end{defi}

\noindent The basepoint in any based set is denoted by $0$.
\smallskip

The category $\Ga$ has a \emph{spine} (minimal full subcategory equivalent to the whole) with objects $[n] = \{0,1,2,\dots,n\}$ for $n\ge 0$.  Likewise, $\Om$ has a spine with objects $\underline{n} = \{1,2,\dots,n\}$ for $n \ge 0$.  In appropriate situations,  we frequently replace an indexing category by its spine and denote it by the same symbol $\Ga$ or $\Om$.

\subsection{Modules} 

Let $K$ be any commutative ground ring.  In our applications $K$ will usually be the ring $\zed$ of integers.

The functors from $\Ga$ to $K$-modules form an abelian category \KGam\ which we shall simply call the category of left $\Ga$-modules (or left $K\Ga$-modules, if we need to make the ring $K$ explicit).  The cofunctors from $\Ga$ to $K$-modules form the category \mKGa\ of right $\Ga$-modules.  

(Alternatively, these module categories can be defined by first forming the ``ring with many objects'' $K\Ga$,  which is the category whose objects are those of $\Ga$ and whose morphism sets are the free $K$-modules $K\!\Hom_\Ga(S,T)$.)

The categories \KOmm\ and \mKOm\ of left and right $\Om$-modules are defined in a way completely analogous to the above.

Replacing the categories $\Ga$ and $\Om$ by their spines changes these module categories only by an equivalence, and therefore for homological purposes makes no difference at all.

\subsection{The Pirashvili-Dold-Kan Theorem for $\Ga$-modules\label{sec: DK}}

We now introduce the cross-effect functors of Eilenberg and Mac Lane, which convert  $\Ga$-modules into $\Om$-modules \cite{pir:dold,pir:hodge}.

Let $A$ be a finite set.  We denote by $A^+$ the based set $A \sqcup \{0\}$, and for each element $a \in A$ we denote by $A^+/\{a\}$ the quotient set obtained by identifying the single element $a$ to the basepoint. Let $r_a: A^+ \to A^+/\{a\}$ be the projection map.

\begin{defi} Let $L$ be a left $\Ga$-module.  The \emph{cross-effect} $\cro(L)$ is the left $\Om$-module which is given on objects by 
$$\begin{aligned} \cro(L)(A) \quad 
&= \quad \bigcap_{a \in A}\; \ker  \left( L(r_a): L(A^+) \to L(A^+/\{a\}) \right)     \\
&= \quad\ker \left(\,\prod_{a \in A}L(r_a): L(A^+) \to  \prod_{a \in A}L(A^+/\{a\})\right) \; ;  \end{aligned}$$
 its value on a morphism $\phi:A \to B$ of $\Om$ is the restriction of \hbox{$L (\phi^+):A^+ \to B^+$}.

If $M$ is a right $\Ga$-module, we define the right $\Om$-module $\cro(M)$ on objects by the dual formula
$$\cro(M)(A) \quad = \quad \coker \left(\,\bigoplus_{a \in A}M(r_a): \bigoplus_{a\in A} M(A^+/\{a\}) \to M(A^+) \right) \; ;$$ 
for a morphism $\phi:A \to B$ of $\Om$ we define $\cro(\phi)$ to be the map of cokernels 
\hbox{$\cro(M)(A) \leftarrow \cro(M)(B)$} induced by $M(A^+) \from M(B^+)$.
\end{defi}
The surjectivity of $\phi$ ensures that $\cro(L)(\phi)$ and $\cro(M)(\phi)$ are well defined \cite{pir:dold}.
\smallskip

These functors give a Morita equivalence between the categories \KGam\ and \KOmm.  A similar result holds for right modules. 

\begin{prop}\emph{(Pirashvili \cite{pir:hodge,pir:dold}} The cross-effect functors 
$$ \cro:{\text \KGam} \quad \stackrel{\sim}{\lra} \quad {\text \KOmm} $$ 
$$ \cro:{\text \mKGa}  \quad \stackrel{\sim}{\lra} \quad {\text \mKOm} $$
are Morita equivalences of abelian categories. \label{prop: DK} $\Box$ \end{prop}

The inverse Morita equivalences $$\crsh \colon {\text \KOmm} \quad \stackrel{\sim}{\lra} \quad {\text \KGam}$$ 
$$\crsh \colon {\text \mKOm} \quad \stackrel{\sim}{\lra} \quad {\text \mKGa}$$
can be described explicitly.  For a left $\Om$-module $L$ there is an isomorphism 
$$\crsh(L)(A^+) \quad \approx \quad \bigoplus_{X \subset A} L(X) \; . $$
When $\psi\colon A^+ \to B^+$ is a map in $\Ga$, its image $\crsh(L)(\psi)$ maps the component labelled $X$ by $L(\psi)$ to the component labelled $\psi(X)$ if $0 \notin \psi(X)$, and to zero otherwise.  The analogous formula holds for right modules \cite{pir:dold}.
\medskip

In this direct sum decomposition for $\crsh(L)(A^+)$ we may take the component corresponding to 
$X = A$.  This gives a natural embedding of $\Om$-modules $$\theta(L)\colon L \to \crsh(L)$$ where the $\Ga$-module $\crsh(L)$ is regarded as an $\Om$-module by restriction along the functor $A \mapsto A^+$.   For a right $\Om$-module $M$ there is likewise a natural projection of $\Om$-modules 
$$\theta(M)\colon \crsh(M) \to M\;.$$

\subsection{Stable homotopy and cohomotopy of $\Ga$-modules\label{sec:pi^*}}

Let $T$ be a right $\Gamma$-module.  The cofunctor $T$ converts any based simplicial set $X$ into a  cosimplicial module $TX$: we denote the cohomology of the associated cochain complex by 
$\pi^*(TX)$.  When $X$ is a sphere, the cohomotopy group $\pi^{n+r}(TS^r)$ is independent of $r$ for $r \ge n+2$ and the common value is called the $n$th {\it stable cohomotopy group} $\pi^nT$ of $T$.

If $F$ is a left $\Ga$-module, then we write $\pi_nF$ for the $n$th \emph{stable homotopy group} of $F$. This is the value of $\pi_{n+r}(FS^r)$ for $r \ge n+2$.

One particularly important right $\Gamma$-module is the linear functor $t$, whose value at a set $S$ is the module of based maps from $S$ to $K$ (where $K$ has basepoint~0).  The right $\Om$-module corresponding to $t$ under the Pirashvili-Dold-Kan isomorphism of \ref{sec: DK} is $\varpi$, where
$$ \varpi(\underline n) = \begin{cases} K& \text{if $n=1$}, \\
                                                  0& \text{if $n \ne 1$}.   \end{cases} $$ 
In topological terms, $t$ is the based $K$-cochain functor.

\begin{thm} \label{thm: cohtpy}  Let $T$ be a right $\Ga$-module,.  Then there are natural isomorphisms
$$ \begin{aligned} \pi_*F \; \approx \; \Tor^\Ga_*(t,F) \; &\approx \; \Tor^\Om_*(\varpi, \cro(F)) \\
  \pi^*T \; \approx \; \Ext_\Ga^*(t,T) \; &\approx \; \Ext_\Om^*(\varpi, \cro(T))  . \end{aligned} $$ \end{thm}

{\it Proof.}  We refer to \cite{pir:dold, pir:hodge} for the proof of this basic result.~$\Box$

\subsection{Deconstructing  the bar construction\label{sec: decon}}

After Theorem \ref{thm: cohtpy} we can calculate stable homotopy $\pi_*F$ and cohomotopy $\pi^*T$ from the bar construction over $\Om$.  Let us write $\til F= \cro(F)$, $\til T=\cro(T)$ for typical left and right $\Om$-modules, and $N\Om$ for the nerve of the category $\Om$.

The module of $q$-chains $\mcl{B}_q(\varpi,\Om,\til F)$  of the bar construction is the direct sum of modules 
 $$\varpi(S_0) \otimes KN\Om_q(S_q,S_0) \otimes \til F(S_q)$$
taken over all $q$-simplices $[f_1|f_2|\dots|f_q]$ of the nerve $N\Om$, where
$$S_0 \stk{f_1}{\lla} S_1 \stk{f_2}{\lla} \dots \stk{f_{q}}{\lla} S_q \quad .$$
Since $\varpi(S_0) = 0$ unless $|S_0| = 1$, and all one-element sets are uniquely isomorphic, we  need only retain terms where $S_0 = \ulin{1}$.  Normalizing by omitting degenerate simplices, we may also assume $|S_1| > 1$.  The face homomorphisms are defined as usual:  $\del_0 = \varpi(f_1) \otimes 1$, whilst  $\del_1, \dots \del_{q-1}$ compose adjacent morphisms in $[f_1|f_2|\dots|f_q]$, and 
$\del_q = 1 \otimes \til F(f_q)$. 

 We can organise $\mcl{B}_*(\varpi,\Om,\til F)$ into a bicomplex as follows.  If $f_{q-t}$ is the last morphism in the sequence $f_1, f_2, \dots , f_q$ which is not an isomorphism (so $t\ge 0$), then we assign the bidegree $(t,q-t)$ to the chain $x \otimes [f_1|f_2|\dots|f_q] \otimes \kappa$.  Then 
$[f_{q-t+1}| \dots|f_{q-1}|f_{q}]$ is a $t$-simplex of the nerve of the subcategory $\Iso(\Om)$ of isomorphisms.  Also $[f_1|f_2|\dots|f_{q-t}]$ is a $(q-t)$-simplex of $N\Om(S_{q-t},\ulin{1})$, and actually a $(q-t-2)$-simplex of $N(S_{q-t}/\!/\Om/\!/\ulin{1})$, which space is isomorphic to 
$N(S_q/\!/\Om/\!/\ulin{1})$ through right-composition with the isomorphism $f_{q-t-1}\dots f_{q-1}f_q$. On such a $(t,q-t)$-chain we define bicomplex differentials by
$$\del' = \sum_{i=q-t}^q (-1)^i \del_i \qquad \text{and} \qquad \del'' = \sum_{i=0}^{q-t-1} (-1)^i \del_i\quad.$$
Then $\del'$ and $\del''$ have degrees $(-1,0)$ and $(0,-1)$ respectively; they anticommute, both square to zero, and their sum is the standard differential.  Thus we have a bicomplex, and therefore a spectral sequence 
$$E^2_{t,s} \approx H'_tH''_s \Longrightarrow H_{t+s} \approx \pi_{t+s}F\quad .$$
Now we calculate the $E^1$ and $E^2$ terms.  To simplify notation, we cut down the category $\Om$ to the spine with objects $\ulin{n} =\{1,2,\dots,n\}$, $n \ge 0$.  This does not prejudice the calculation.  The $E^0$ term is the chain bicomplex.  Since the category $\Iso(\Om)$ now reduces to the disjoint union of the symmetric groups $\Sigma_n$, and $ \varpi(\ulin{1}) \approx K$, we have by the above 
$$E^0_{t,s} \quad = \quad \bigoplus_{n\ge 2}\til C_{s-2}(N(\ulin{n}/\!/\Om/\!/\ulin{1}) \otimes C_t(\Sigma_n; \til F(\ulin{n})) \quad .$$

The $E^1$ term is the $\del''$-homology of $E^0$.  Since $\del_0 =0$ under our normalization, the differential $\del''$ is exactly the differential in the augmented chain group 
$\til C_{s-2}(N(\ulin{n}/\!/\Om/\!/\ulin{1})$.  The homology of this is, by Lemma~\ref{lemma: slicetree} and Proposition~\ref{prop: partbdry}, non-trivial only when $s-2 = n-3$, and then it is isomorphic as a right $\Sigma _n$-module to the dual superlie representation $\mcl{S}^*_n$ of Definition~\ref{def: super}.  This means that $(E^1,d^1)$ is the two-sided bar construction over $\Iso(\Om)$
$$ E^1_{s,t} \; \approx \; \mcl{B}_t(\mcl{S}^*_{s+1}, \Sigma_{s+1}, \til F(\ulin{s+1}))$$
and $E^2_{s,t} \approx \Tor_t^{\Sgm_{s+1}}(\mcl{S}_{s+1},\til F(\ulin{s+1}))$.  We have proved the first statement in the following theorem.  The second part is proved analogously.

\begin{thm} \label{thm: bicomplex}  Let $F$ be any left $\Ga$-module, and $\til F = \cro(F)$.  Then there is a convergent spectral sequence 
$$E^2_{s,t} \approx \Tor_t^{\Sgm_{s+1}}(\mcl{S}^*_{s+1},\til F(\ulin{s+1})) \Longrightarrow \pi_{s+t}F$$ 
where $\mcl{S}^*_{s+1}$ is the dual superlie representation.  If $T$ is a right $\Ga$-module, and $\til T = \cro(T)$, then there is a convergent spectral sequence 
$$E_2^{s,t} \approx \Ext^t_{\Sgm_{s+1}}(\mcl{S}^*_{s+1},\til T(\ulin{s+1})) \Longrightarrow \pi^{s+t}T\quad.\quad\Box $$ 
\end{thm}
\medskip

\begin{cor} \label{cor: linecor} If $\til F(\ulin{i}) = 0$ for all $i$ except for a single value $i = n$ then 
$$\pi_k(F) \approx \Tor_{k-n+1}^{\Sgm_{n}}(\mcl{S}^*_{n},\til F(\ulin{n})) \quad \text{for all $k$.}$$
If $\til T(\ulin{i}) = 0$ for all $i$ except for $i = n$, then $$\pi^kT \approx \Ext^{k-n+1}_{\Sgm_{n}}(\mcl{S}^*_{n},\til F(\ulin{n}))\quad \text{for all $k$}.  $$ \end{cor}

{\it Proof.} The spectral sequence of Theorem~\ref{thm: bicomplex} collapses to one line. $\Box$
\bigskip

We note that all simple $\Om$-modules are of the type described in Corollary \ref{cor: linecor} \cite{pir:hodge}.
\smallskip

The following further corollaries will be useful in the proof of Theorem \ref{thm: Gacoho} in the next section.  We consider the situation where the only non-zero value~$\til F(\ulin{n})$ of the $\Om$-module $\til F$ is an induced module $\Ind_{\Sigma_{n-1}}^{\Sigma_n}\!M$ for some $\Sigma_{n-1}$-module~$M$.

\begin{cor} \label{cor: triv} If $M$ a $\Sigma_{n-1}$-module, and $F$ is a left $\Ga$-module with the property 
$$\cro(F)(\ulin{i})  =  \til F(\ulin{i}) = \begin{cases}   \Ind_{\Sigma_{n-1}}^{\Sigma_n}\!M &  \text{(as a $\Sigma_n$-module) when $i=n$}, \\
0  & \text{    when $i \ne n$} \end{cases} $$
then $\pi_{n-1}(F) \approx M$, and $\pi_k(F) \approx 0$ for all $k \ne n-1$. \\
\\
If $M$ is a right $\Sigma_{n-1}$-module, and $T$ is a right $\Ga$-module with the property 
$$\cro(T)(\ulin{i})  =  \til T(\ulin{i})  =  \begin{cases}   \Ind_{\Sigma_{n-1}}^{\Sigma_n}\!M &  \text{(as a $\Sigma_n$-module) when $i=n$}, \\
0  & \text{    when $i \ne n$} \end{cases} $$
then $\pi^{n-1}(T) \approx M$, and $\pi^k(T) \approx 0$ for all $k \ne n-1$.
 \end{cor}

{\it Proof.}  In the case of left modules, Corollary \ref{cor: linecor} and a standard induction formula give 
$$\begin{aligned}         \pi_k(F) \quad \approx \quad 
                    &\Tor_{k-n+1}^{\Sgm_{n}}(\mcl{S}^*_{n},\til F(\ulin{n})) \\
        \quad  \approx \quad 
&\Tor_{k-n+1}^{\Sgm_{n}}(\mcl{S}^*_{n}, \Ind_{\Sigma_{n-1}}^{\Sigma_n}\!M) \\
        \quad  \approx \quad   
&\Tor_{k-n+1}^{\Sgm_{n-1}}(\Res_{\Sigma_{n-1}}^{\Sigma_n}\!\mcl{S}^*_{n}, M) \\
       \quad  \approx \quad 
&\Tor_{k-n+1}^{\Sgm_{n-1}}(K\Sigma_{n-1}, M) \\
       \quad   \approx \quad   & \begin{cases} M &\text{if $k = n-1$} \\
                                            0  &\text{if $k \ne n-1$}  \end{cases} \end{aligned} $$
since $\mcl{S}^*_n$ restricts to the regular representation of $\Sigma_{n-1}$ by Proposition~\ref{prop: partbdry}.
The case of right modules is completely analogous. $\Box$

\subsection{\label{sec: RW}The complexes $\Cdnga_*(F), \Cupga^*(T)$ }

The stable homotopy and cohomotopy groups of $\Ga$-modules will arise geometrically in our operadic obstruction theory as the $\Tor$ and $\Ext$ groups of Theorem~\ref{thm: cohtpy}.  The geometry of operads does not directly give the bar construction, however, but rather a significant variant of this, namely a reduced version $\Rupga^*(T)$ of the Robinson-Whitehouse complex $\Cupga^*(T)$ of \cite{rob-whi2,pir-richt}.   We therefore require the main theorem of  \cite{pir-richt}, which shows that this complex does indeed represent stable homotopy and cohomotopy.

We begin with the case of left modules and homotopy, for which we adapt the definition given in \cite{pir-richt}.  For any left $\Ga$-module $F$ we define a chain complex $\Cdnga_*(F)$ as follows.  We regard $\Om$ as embedded in $\Ga$ by the basepoint-adjunction functor $S \mapsto S^+$. Just as in the bar construction of \ref{sec: decon}, for each unbased finite set $S$ let $N\Om_q(S,\underline{1})$ be the set of  $q$-simplices 
$[f_1|f_2|\dots|f_q]$ of the nerve of~$\Om$
$$\underline{1} \approx S_0 \stk{f_1}{\lla} S_1 \stk{f_2}{\lla} \dots \stk{f_{q}}{\lla} S_q=S$$
which have first source $S$, and have final target a set with one element which we may take to be $\underline{1}$.  For $n\ge 2$, $N\Om_q(S,\underline{1})$ is therefore just the set of $(q-2)$-simplices of the slice category $S/\Om/\underline{1}$.  Let $KN\Om_q(S,\underline{1})$ be the free $K$-module generated by $N\Om_q(S,\underline{1})$, and let 
$$\Cdnga_q(F) \quad = \quad \bigoplus_SKN\Om_q(S,\underline{1}) \otimes F(S^+) \;.$$
Face operators
$\partial_i\colon \Cdnga_q(F) \to \Cdnga_{q-1}(F)$ are defined as follows.
$$\begin{aligned}
\partial_q([f_1|f_2|\dots|f_q] \otimes x) \quad &= \quad [f_1|f_2|\dots|f_{q-1}] \otimes F(f_q^+)(x) \\
\partial_i([f_1|f_2|\dots|f_q] \otimes x) \quad &= \quad [f_1|\dots|f_if_{i+1}|\dots|f_q] \otimes x
\end{aligned}$$
for $0 < i < q$.  Note that our embedding of $\Om$ in $\Ga$ takes $f_q$ into a morphism $f_q^+$, and  $F(f_q^+)$ is the value of $F$ on this.  To define $\partial_0$, partition the source $S_1$ of $f_1$ into singleton sets~$\{s\}$.  Taking inverse images induces partitions $S_j = \bigsqcup_{s \in S_1} S_j^{(s)}$ of the sets $S_2,\dots,S_q$, where $S_j^{(s)} = (f_2f_3\dots f_j)^{-1}\{s\}$.  Let $f_j^{(s)}\colon S_j^{(s)} \to S_{j-1}^{(s)}$  be the restriction of $f_j$, and let $r_s\colon S_q^+ \to S_q^{(s)+}$ be the pointed map which is the identity on $S_q^{(s)}$ and maps its complement to the basepoint. We set 
$$\partial_0([f_1|f_2|\dots|f_q] \otimes x) \quad 
= \quad \sum_{s\in S_1}[f_2^{(s)}|\dots|f_q^{(s)}] \otimes F(r_s)(x)\; .$$
It is straightforward to verify that, for any left $\Ga$-module $F$, the $\partial_i$ satisfy the usual identities for face operators, giving a simplicial $K$-module.  The associated chain complex, with boundary $\sum (-1)^i\partial_i$, is the Robinson-Whitehouse complex $\Cdnga_*(F)$.  Its homotopy type is unchanged when the category $\Ga$ is replaced by a spine, so we may always do this when computing its homology, which we denote by $H^\Ga_*(F)$.

We may now choose $F$ to be the $\Ga$-bimodule $K\Ga$.  Then $\Cdnga_*(K\Ga)$ is a complex of right $\Ga$-modules.  The definition \ref{sec:pi^*} of the right $\Ga$-module $t$ gives an isomorphism $H^\Ga_0(K\Ga) \approx t$.

We define the Robinson-Whitehouse cochain complex $\Cupga^*(T)$ of a right $\Ga$-module $T$ by $$\Cupga^*(T) \quad \approx \quad \Hom_\Ga(\Cdnga\!(K\Ga), T)\;.$$
We denote its cohomology by $H_\Ga^*(T)$.  The following result contains the main result of~\cite{pir-richt} and its cohomological dual.

\begin{thm} \label{thm: Gacoho}
The complex $\Cdnga_*(K\Ga)$ is a projective resolution of the right $\Ga$-module $t$.
 There are natural isomorphisms
$$\begin{aligned} H^\Ga_*(F) \quad \approx \quad &\Tor^\Ga_*(t, F) \quad \approx \quad \pi_*\!(F) \\
H_\Ga^*(T) \quad \approx \quad &\Ext_\Ga^*(t, T) \quad \approx \quad \pi^*\!(T) \end{aligned} $$
for all left $\Ga$-modules $F$ and right $\Ga$-modules $T$.
 \end{thm}

Following Lemma \ref{lemma: gaspec} below, we shall prove Theorem~\ref{thm: Gacoho} by reducing it to Theorem~\ref{thm: cohtpy}.  
\medskip

\noindent {\bf Notes.}  \begin{enumerate}
\item We observe that the definition of the complexes $\Cdnga_*(F)$ and $\Cupga^*(T)$ uses only the $\Theta$-structure of $F$ or $T$, where $\Theta \subset \Ga$ is the subcategory of morphisms $f \colon S \to S'$ in $\Ga$ which are surjections of based sets.  (Thus $\Theta$ is generated by the embedded copy of $\Om$ together with the morphisms $r_\alpha$ of \ref{sec: DK}.)  In fact, the complexes are defined for every $\Theta$-module.  We shall require this fact just once, in the proof of the following lemma.
\item The  complexes $\Cdnga_*(F)$ and $\Cupga^*(T)$ are exact functors of $F$ and $T$ respectively, and so $H^\Ga_*(F)$ and $H_\Ga^*(T)$ are connected sequences of (co)homological functors.
\end{enumerate}

\begin{lemma} \label{lemma: gaspec} If $F$ is a left $\Ga$-module there is a natural spectral sequence 
$$E^1_{p,q} \approx \Tor_q^{\Sgm_{p+1}}(\mcl{S}^*_{p+1}, F[p+1]) 
\Longrightarrow H^\Ga_{p+q}(F) \;.$$
If $T$ is a right $\Ga$-module there is a natural spectral sequence 
$$E_1^{p,q} \approx \Ext^q_{\Sgm_{p+1}}(\mcl{S}^*_{p+1}, T[p+1]) 
\Longrightarrow H_\Ga^{p+q}(T) \;.$$
\end{lemma}

{\it Proof.}  We take the case of a left module $F$.  We can filter $F$ by $\Theta$-submodules $F^n$, where 
$$F^n[j] \quad = \quad \begin{cases}  F[j] & \text{if $j \le n$} \\
                                              0 &  \text{if $j > n$} \end{cases} $$
and form the quotient $\Theta$-modules $F^j/F^i$ for $j \le i$.  In view of the Notes above, this gives a convergent spectral sequence
$$E^1_{p,q} \quad \approx H^\Ga_{p+q}(F^{p+1}/F^p)  \quad
\Longrightarrow \quad H^\Ga_{p+q}(F)\;.$$
Since $F^{p+1}/F^p$ is concentrated in the single degree $p+1$, the chain complex $\Cdnga_*(F^{p+1}/F^p)$ has $\del_0 = 0$ for nondegenerate chains, and is isomorphic with the bar resolution $\mcl{B}(\varpi,\Om,F^{p+1}/F^p)$ of \ref{sec: decon}. 
We can therefore read off the $E^1$ term from Corollary \ref{cor: linecor}:
$$\begin{aligned} E^1_{p,q} \quad &\approx \quad H^\Ga_{p+q}(F^{p+1}/F^p) \\
                                & \approx \quad \Tor_{p+q}^\Om(\varpi, (F^{p+1}/F^p)) \\
                                & \approx \quad \Tor_{q}^{\Sigma_{p+1}}(\mcl{S}^*_{p+1}, F[p+1])\;.
\end{aligned} $$
The differential $d^1_{p,q} \colon E^1_{p,q} \to E^1_{p-1,q}$ is a sum over principal faces of the operad space $\mcl{T}_{p+1}$.  In algebraic terms it is a composite
$$\Tor_{q}^{\Sigma_{p+1}}(\mcl{S}^*_{p+1}, F[p+1]) \lra
\Tor_{q}^{\Sigma_{p-1} \times \Sigma_2}(\mcl{S}^*_{p+1}, F[p+1])  \lra
\Tor_{q}^{\Sigma_{p}}(\mcl{S}^*_{p}, F[p]) $$
where the first homomorphism is the homology transfer.   The second homomorphism combines a boundary map in the homology of $\mcl{S}_{p+1}$, and  $F[p+1] \to F[p]$ induced by a surjection $[p+1] \to [p]$.

The spectral sequence for a right module $T$ is established similarly.~$\Box$
\bigskip

{\bf Proof of Theorem \ref{thm: Gacoho}.}  We first prove that for left $\Ga$-modules~$F$ there is a natural isomorphism  $\theta(F)\colon \Tor^\Om_*(\varpi, \til F) \to H^\Ga_*(F)$, where again $\til F = \cro(F)$ is the $\Om$-module corresponding to $F$. After Theorem \ref{thm: cohtpy}, this will establish one main statement of Theorem \ref{thm: Gacoho}. 

As  usual we embed $\Om$ into $\Ga$ by adding a basepoint.  There is an inclusion $\mcl{B}(\varpi,\Om,\til F) \lra \Cdnga(F)$ induced by the inclusion of $\til F$ into $F$.  This commutes with the face operator $\del_0$ because $\til F$ lies in the kernel of each projector $r_\alpha$ of (\ref{sec: DK}); it commutes with other face operators for trivial reasons.  Therefore it is a chain map, and induces homomorphisms of connected sequences of homology functors $\theta(F)\colon \Tor^\Om_*(\varpi, \til F) \to H^\Ga_*(F)$.  We show that $\theta(F)$ is an isomorphism.

\smallskip \emph{Case 1.  The module $\til F$ is concentrated in a single degree $n$.} 

In this case, $\mcl{B}(\varpi,\Om,\til F)$ maps isomorphically onto $\Cdnga(F^n)$, where $F^n$ is the first non-trivial truncation of $F$ (as a $\Theta$-module, in Lemma~\ref{lemma: gaspec}).  Therefore it is enough to show that $H^\Ga(F/F^n)$ is zero.  We use the spectral sequence of Lemma~\ref{lemma: gaspec}) for the module $F/F^n$.  If $F[n]$ is the $\Sigma_n$-module $M$, then the formula for the inverse Morita equivalence (following Proposition~\ref{prop: DK}) shows that 
$F[n+k] \approx \Ind_{\Sigma n}^{\Sigma n+k}\!M$.  Thus all the modules in $F/F^n$ are induced representations, and it follows exactly as in Corollary~\ref{cor: triv} that all the higher $\Tor$-groups vanish in our spectral sequence 
$$E^1_{p,q} \approx \Tor_q^{\Sigma_{p+1}}(\mcl{S}^*_{p+1}, (F/F^n)[p+1]) 
\Longrightarrow H^\Ga_{p+q}(F/F^n) $$
which therefore collapses to the edge $q=0$.  Further, this edge has groups non-zero only when $p \ge n$.   These are 
$$\begin{aligned}
E^1_{p,0} \quad \approx& \quad
\mcl{S}^*_{p+1} \otimes_{\Sigma_{p+1}} \Ind_{\Sigma_n}^{\Sigma_{p+1}}\!M \\
                  \approx& \quad \Ind_{\Sigma_n}^{\Sigma_{p}}\!M \end{aligned}\;.$$
This $(E^1,d^1)$ is a cone.  It has a contracting homotopy, given by inclusion homomorphisms $\Ind_{\Sigma_n}^{\Sigma_{p}}\!M \to \Ind_{\Sigma_n}^{\Sigma_{p+1}}\!M \;.$  Therefore 
$H^\Ga(F/F^n)$ is zero, and $\theta(F)$ is an isomorphism. 

\smallskip \emph{Case 2.  The $\Om$-module $\til F$ has finite length.}
Since simple $\Om$-modules are concentrated in a single degree, this follows from Case~1 by exactness, the five-lemma and finite induction.

\smallskip \emph{Case 3.  The module $\til F$ is a general left $\Om$-module.}
This follows from Case 2 by taking an inductive limit.

\smallskip It follows as usual that $\Cdnga(K\Ga)$ is a projective resolution of $\varpi$.

\medskip The cohomological statement in Theorem~\ref{thm: Gacoho} is equivalent to 
$$H_\Ga^*(T) \quad \approx \quad \Ext_\Om^*(\varpi, \til T)$$ naturally in the right $\Ga$-module $T$.  The proof is dual to that above.  There is no extra problem in the final stage, as mapping from the locally finite object $\Cdnga(\Ga)$ commutes with taking a filtered inductive limit in the target variable $T$.

This completes the proof of Theorem~\ref{thm: Gacoho}.  $\Box$
\medskip

The original proof of Theorem~\ref{thm: Gacoho}, due to Pirashvili and Richter~\cite{pir-richt}, relates the complex $\Cdnga_*(F)$ directly to stable homotopy.

\subsection{\label{sec: redcx}The normalized complexes $\Rdnga_*(F), \Rupga^*(T)$}
\phantom{B}

The complex $\Cdnga_*(F)$ in \S\ref{sec: RW} is that of \cite{pir-richt,rob-whi2}. It does not exactly match the obstruction theory which we develop in \S\ref{sec: obstructions}, because it has no diagonal filtration  by subcomplexes (see Definition~\ref{defi: filt} below).  We therefore introduce a normalized subcomplex $\Rdnga_*(F)$.  Whereas 
$$\Cdnga_q(F) \quad = \quad \bigoplus_SKN\Om_q(S,\underline{1}) \otimes F(S^+)\;, $$
in $\Rdnga_*(F)$ we require only those indexing simplices $[f_1|f_2|\dots|f_q]$ of $N\Om_q(S,\ulin{1})$ in which isomorphisms precede strict epimorphisms.

 As the object $\ulin{1}$ is terminal in $\Om$, we have $N(S/\Om) \approx N(S/\Om/\ulin{1})$.  We now define the width of a $(q-1)$-simplex $[f_1|f_2|\dots|f_q]\eps$ of $N(S/\Om)$. 

\begin{defi} We say the simplex 
$$\underline{1} \stk{f_1}{\lla} S_1 \stk{f_2}{\lla} \dots \stk{f_{q}}{\lla} S_q = S$$
has \emph{width $p$} if all of $f_q, f_{q-1},\dots, f_{q-p+1}$ are isomorphisms in $\Om$, but $f_{q-p}$ is not an isomorphism.\label{defi: width} 

This gives a filtration on $KN\Om(S,\ulin{1})$: we define $\mcl{W}^pKN\Om(S,\ulin{1})$ to be the submodule generated by all simplices of width less than or equal to~$p$. \end{defi}

\begin{defi}\label{def: rednerv} The \emph{reduced nerve} $\til N(S/\Om)$ is the subspace of $N(S/\Om)$ consisting of those simplices $[f_1|f_2|\dots|f_q]$ as above in which {isomorphisms always precede non-isomorphisms}: that is, if the simplex has width~$p$, then $f_{q-p}, f_{q-p-1},\dots, f_2, f_1$ are all strict epimorphisms. 

We use the notation $\til N\Om_q(S,\ulin{1})$ for the set of such simplices.

We define $\Rdnga_q(F)$  to be the submodule  $\bigoplus_SK\til N\Om_q(S,\underline{1}) \otimes F(S^+)$ of~$\Cdnga_q(F)$. \end{defi} 

The face operators in $\Cdnga_*(F)$ respect this submodule, so that $\Rdnga_*(F)$ is a chain subcomplex.

\begin{prop} The inclusion $\Rdnga_*(F) \to \Cdnga_*(F)$ induces isomorphisms of homology for every left $\Ga$-module $F$.  In particular, $\Rdnga_*(K\Ga)$ is a projective resolution of the right $\Ga$-module $t$. \end{prop}

{\it Proof.} Here we filter both complexes $\Rdnga_*(F)$ and $\Cdnga_*(F)$ by height; that is, by the cardinality of the initial set $S$ in the simplex of the nerve of $\Om$.  This gives strongly convergent spectral sequences for the homology of each complex.  In the filtration quotients, the anomalous face operator $\del_0$ of \ref{sec: RW} reduces to the usual $\del_0$ of the bar construction.  Thus again we can apply the bicomplex deconstruction of \ref{sec: decon} and find that for each complex the $E^2$-term is (compare Theorem~\ref{thm: bicomplex})
$$E^2_{p,q} \; \approx \; \Tor_q^{\Sigma_{p+1}}(\mcl{S}^*_{p+1}, \til F(\ulin{p+1}))\;.$$
The difference is that the same homology $\mcl{S}^*_{p+1}$ of $N(\ulin{p+1}/\!/\Om/\!/\ulin{1})$ is calculated in $\Rdnga_*(F)$ from the normalized chain complex, and in $\Cdnga_*(F)$  from the full chain complex.  The inclusion map induces the indicated isomorphism on $E^2$, and therefore an isomorphism on $E^\infty$ and on the abutment. $\Box$

Theorem~\ref{thm: Gacoho}   therefore has the following corollary, in which 
$\Rupga^*(T) = \Hom_\Ga(\Rdnga_*(K\Ga), T)$.

\begin{cor}  There are natural isomorphisms
$$\begin{aligned} H_*(\Rdnga_*(F)) \quad \approx \quad &\Tor^\Ga_*(t, F) \quad \approx \quad \pi_*\!(F) \\
H^*(\Rupga^*(T)) \quad \approx \quad &\Ext_\Ga^*(t, T) \quad \approx \quad \pi^*\!(T) \end{aligned} $$
for all left $\Ga$-modules $F$ and right $\Ga$-modules $T$.
 \end{cor}

We now use the width filtration of Definition~\ref{defi: width} to construct  another filtration of $\Rdnga_*(F)$ by diagonal subcomplexes. 
\begin{defi} \label{defi: filt}The $n$th stage in the diagonal filtration of $\Rdnga_*(F)$ is the submodule 
$$\Delta^n\Rdnga_q(F) \quad = \quad \bigoplus_{p + |S| \le n}\mcl{W}^pK\til N\Om_q(S,\underline{1}) \otimes F(S^+) \;.$$  
The dual complex $\Rupga^*(T)$ of a right $\Ga$-module $T$ has a dual filtration~in which 
$\Delta_n\Rupga^q (T)$ consists of those cochains $\phi\colon \Rdnga_q(K\Ga) \to T$ which annihilate $\Delta^{n-1}\Rdnga_q(K\Ga)$.
\end{defi}
\medskip

\begin{lemma} \label{lemma: filt} The diagonal filtration of $\Rdnga_*(F)$ is a filtration by subcomplexes. \end{lemma}
{\it Proof.} We show that $\Delta^n\Rdnga_*(F)$ is closed under every face operator  of $\Rdnga_*(F)$, whatever the left $\Ga$-module $F$.  This is clearly true for the operators $\del_i$ when $i>0$, since these are non-increasing on both width $p$ and height $|S|$.  There remains $\del_0$.  We recall (see \S\ref{sec: RW}) that $\del_0$ maps a generator of width $p$ 
$$\underline{1} \stk{f_1}{\lla} S_1 \stk{f_2}{\lla} \dots \stk{f_{v}}{\lla} S_v= S$$ to a sum of terms obtained by taking inverse images of the various elements of $S_1$.  Any of these terms, such as that corresponding to a component $S^{(s)}$ say of $S$, may have width $p+i > p$, if $i$ strict epimorphisms restrict to isomorphisms on the $s$-component.  But then these would restrict to $i$ strict epimorphisms on the union of the other components, which implies the increased-width term has height at most $|S| - i - 1$.~$\Box$

The analogue for the full complex $\Cdnga_*(F)$ of the above lemma is \emph{false.} That is the reason for using the reduced version.  The principal property of the above algebraic diagonal filtration is the following.
\medskip

\begin{prop} \label{prop: collapse}Let $F$ be any left $\Ga$-module.  The spectral sequence for  $\pi_*(F)$ arising from the diagonal filtration $\Delta^*\Rdnga_*(F)$ collapses from the $E^2$-term onwards, and the filtration on the abutment is trivial in each degree.

For a right $\Ga$-module $T$ the corresponding statements hold: the spectral sequence converging to $\pi^*(T)$ arising from the filtration $\Delta_*\Rupga^*(T)$ is likewise trivial.
\end{prop}
 
{\it Proof.} For any left $\Ga$-module $F$ the $E^0$ term splits as a sum over values of the height $k$.  If we reduce $\Om$ to its spine,  we can write it in terms of the bar construction over $\Sigma_k$ as
$$\begin{aligned}
E^0_{p,q} \quad \approx& \quad \Delta^p(\Rdnga_{p+q}(F))/\Delta^{p-1}(\Rdnga_{p+q}(F))   \\
    \approx& \quad \bigoplus_{k=1}^p \mcl{B}_{p-k}
\biggl(\til C_{k+q-2}\bigl(
N(\ulin{k}/\!/\Om/\!/\ulin{1})\bigr), \,\Sigma_k,\,F[k]\biggr) \end{aligned}$$
and $d^0_{p,q}\colon E^0_{p,q} \to E^0_{p,q-1}$ is the differential in $\til C_*(N(\ulin{k}/\!/\Om/\!/\ulin{1}))$.  By Lemma~\ref{lemma: slicetree} and Proposition~\ref{prop: partbdry} this has just one non-trivial homology group, namely the module $\mcl{S}_k^*$ in degree $k-3$.  Therefore the non-trivial $E^1$ groups are confined to the single line $q = -1$, where we have 
$$E^1_{p, -1} \quad \approx \quad \bigoplus_{k=1}^p \mcl{B}_{p-k}(\mcl{S}_k^*, \Sigma_k, F[k])\;.$$  Hence the only possible non-zero differential is $d^1$: all differentials from $d^2$ onwards must be zero, and the spectral sequence collapses, with no extension problems.

The analogous argument works in the spectral sequence for the cohomotopy of a right $\Ga$-module $T$. $\Box$

\begin{cor}\label{cor: finalmess} If $F$ is a left $\Ga$-module, then the stable homotopy $\pi_n(F)$ is the homology 
$\;\ker d^1_{n+1,-1}/\im d^1_{n+2,-1}$ of the complex 
$$\begin{aligned} \lra \bigoplus_{k=1}^{n+2} \mcl{B}_{n+2-k}(\mcl{S}_k^*, \Sigma_k, F[k]) \stk{d^1_{n+2,-1}}{\lra} \bigoplus_{k=1}^{n+1}& \mcl{B}_{n+1-k}(\mcl{S}_k^*, \Sigma_k, F[k]) \;\;\lra \\
\stk{d^1_{n+1,-1}}{\lra}&\bigoplus_{k=1}^{n} \mcl{B}_{n-k}(\mcl{S}_k^*, \Sigma_k, F[k]) \lra \end{aligned}$$
where $d^1_{**}$ is induced by the differential in $\Rdnga_*(F)$. \\

If $L$ is any right $\Ga$-module, then the stable cohomotopy $\pi^n(L)$ is the cohomology 
$\;\ker d_1^{n+2,-1}/\im d_1^{n+1,-1}$  of the complex 
$$\begin{aligned} \lla \bigoplus_{k=1}^{n+2} \mcl{B}^{n+2-k}(\mcl{S}_k^*, \Sigma_k, L[k]) \;\stk{d_1^{n+2,-1}}{\lla}\; \bigoplus_{k=1}^{n+1}& \mcl{B}^{n+1-k}(\mcl{S}_k^*, \Sigma_k, L[k]) \;\;\lla \\
\stk{d_1^{n+1,-1}}{\lla}\;&\bigoplus_{k=1}^{n} \mcl{B}^{n-k}(\mcl{S}_k^*, \Sigma_k, L[k]) \lla \end{aligned}$$ where $d_1^{**}$ is induced by the differential in the complex $\Rupga^*(L)$. $\Box$
\end{cor}

\noindent We recall the notation used in the above:  $\mcl{B}_{r}(\mcl{S}_k^*, \Sigma_k, F[k])$ denotes the bar resolution module 
$\mcl{S}_k^* \otimes (\zed\Sigma_k)^{\otimes r} \otimes F[k]$, whilst $\mcl{B}^{r}(\mcl{S}_k^*, \Sigma_k, L[k])$ denotes $\Hom(\mcl{S}_k^* \otimes (\zed\Sigma_k)^{\otimes r}, \,L[k])$.  \medskip

The chain complex in Corollary~\ref{cor: finalmess} is equivalent to the total $\Xi$-complex of~\cite{rob:inv}, and the above result is equivalent to Corollary~3.7 of that paper.

\section{Obstruction theory for operads \label{sec: obstructions}}

\subsection{\label{sec: piop}A $\Ga$-module structure on the homotopy of an operad}

From now on (unless otherwise stated) our commutative ground ring $K$ will be either the ring $Z$ of integers, or a suitable localization or completion of $Z$.  In the complete case all modules are assumed formally complete.

We now introduce the basic product structure in an operad $\mcl{C}$ which we shall hope to refine to an $E_\infty$ structure.  The operad $\mcl{C}$ is required to satisfy a fibrancy condition which will be needed in \S\ref{sec: obstructions}.

\begin{defi}  Let $\mcl{C}$ be an Kan operad; that is, each $\mcl{C}_S$ is a Kan complex.  
\label{def: mon}
An \emph{h-monoid} in $\mcl{C}$ is a pair $(\mu,\eta)$ where 
\begin{enumerate} \item \quad $\mu$ is a point in $\mcl{C}_2$, \ \emph{and} 
                           \item \quad $\eta$ is a point $\mcl{C}_0$, \ \emph{such that} 
                           \item \quad $\mu \circ_1 \mu$  and $\mu \circ_2 \mu  \quad$ are in the same path-component of $\:\mcl{C}_3$;
                            \item \quad $\mu$ and $\sigma(\mu)$ are in the same path-component of $\mcl{C}_2$, where $\sigma$ generates~$\Sigma_2$; and
                           \item \quad the element $\iota = \mu \circ_2 \eta\; \in \; \mcl{C}_1$ is a 
\emph{homotopy unit:} any left or right composition with $\iota$ is homotopic to the identity map on any $\mcl{C}_S$.  \end{enumerate} \end{defi}

\noindent This definition contains minor abuses of notation.  In (5) above, for instance, $\iota \in \mcl{C}_1$ has been identified with its canonical image in $\mcl{C}_{\{s\}}$ for any $s\in S$.  In following sections we shall continue to elide isomorphisms by calling them relabellings.

\begin{defi} \label{looperad} We say that $(\mu,\eta)$ is a \emph{loop}  h-monoid in $\mcl{C}$ if each space $\mcl{C}_n$ is homotopy equivalent to a loop space, and  each composition with $\mu$ or $\eta$ is homotopic to a loop map.
\end{defi}

\begin{lemma} \label{lemma: gastr} Let $\mcl{C}$ be an Kan operad, and $(\mu,\eta)$ a loop h-monoid in $\mcl{C}$.

Then $(\mu,\eta)$ induces a right $\Ga$-module structure on $\pi_k\mcl{C}$ for all $k\ge1$.
\end{lemma}

{\it Proof.} Since $\mcl{C}_n$ is a loop space, it is a simple space, and  its homotopy groups are independent of basepoint (up to unique isomorphism).  Its fundamental group is abelian, and any unbased map $S^k \to \mcl{C}_n$ gives a well-defined element of $\pi_k\mcl{C}_n$.

We use the h-monoid $(\mu,\eta)$ to define \hbox{$\phi^*\colon\pi_k\mcl{C}_n \to \pi_k\mcl{C}_m$} for each morphism \mbox{$\phi\colon [m] \to [n]$} in $\Ga$, and to show that $(\phi\psi)^* = \psi^*\phi^*$.  

For an isomorphism $\sigma \in \Sigma_n \subset \Ga$, we define the action $\sigma_*$ on $\pi_k\mcl{C}_n$ to be that induced by the right action of $\Sigma_n$ on $\mcl{C}_n$ in the operad $\mcl{C}$.  Every other morphism $\phi$ in $\Ga$ can be written as $\lambda\circ\omega\circ\tau$, where 
\begin{enumerate} \item[] \quad $\lambda$ is injective
                           \item[] \quad $\omega$ is surjective and satisfies $\omega^{-1}(0) = \{0\}$
                           \item[] \quad $\tau$ is surjective, and bijective away from $\tau^{-1}(0)\,$;
\end{enumerate} and this factorization is unique up to isomorphisms.  Therefore $\Ga$ is generated, modulo isomorphisms, by the following elements, the action of each of which we specify.   \begin{enumerate}
\item The inclusion $\lambda\colon [n] \to [n+1]$.   We define  
$\lambda^*\colon \pi_k\mcl{C}_{n+1} \to \pi_k\mcl{C}_n$ to be the morphism induced by composition 
with $\eta$ in the $(n+1)$st place $$(\frac{\;\;}{\;\;\;}\circ_{n+1} \eta)\colon \mcl{C}_{n+1} \to \mcl{C}_n.$$  

\item The surjection $\omega:[n+1] \to [n]$ which maps $n+1$ to $n$, and is elsewhere the identity map.
We define $\omega^*\colon \pi_k\mcl{C}_{n} \to \pi_k\mcl{C}_{n+1}$ to be the morphism induced by right  composition with $\mu$ in the $n$th place
$$(\frac{\;\;}{\;\;\;}\circ_n \mu)\colon \mcl{C}_{n} \to \mcl{C}_{n+1}$$
the inputs to $\mu$ being relabelled $n$ and $n+1$. 
\item The surjection $\tau\colon [n+1] \to [n]$ which maps $n+1$ to the basepoint $0$, and is elsewhere the identity map.  We define $\tau^*\colon \pi_k\mcl{C}_{n} \to \pi_k\mcl{C}_{n+1}$ to be the morphism induced by left composition with $\mu$
$$(\mu\circ_2 \frac{\;\;}{\;\;\;})\colon \mcl{C}_{n} \to \mcl{C}_{n+1}$$ 
the other input of $\mu$ being relabelled $n+1$.
\end{enumerate}
We assumed that the operad $\mcl{C}$ consists of loop spaces.  The isomorphisms between homotopy groups of different components of a loop space are induced by loop composition.  Therefore the hypothesis that compositions with $\mu$ and $\eta$ are loop maps, implies that each of these induced homomorphisms of $\pi_k$ is independent of basepoint.

In order that these constructions make $\pi_k\mcl{C}$ into a right $\Ga$-module, they must respect the relations among $\lambda$, $\omega$, $\tau$ and the permutations 
in~$\Ga$.  All of these follow from the associativity laws for an operad (Def.\ref{def: operad}, clauses 3,4) together with the relations in an h-monoid (Def.\ref{def: mon}, clauses 3,4,5).   For instance, as $\tau\lambda =1$ in $\Ga$, we need that $\lambda^*\tau^* =1$; but this follows from right-right associativity in $\mcl{C}$ (Def.\ref{def: operad}, clause~4) together with the homotopy unit condition for an h-monoid  (Def.\ref{def: mon}, clause~5).~$\Box$

\subsection{Spaces of  $E_\infty$ structures on an operad\label{sec: modsp}}

We investigate the problem of mapping the universal $E_\infty$ operad $\mcl{T}$ into a general operad~$\mcl{C}$.

\begin{defi} \label{def: formal E} Let $\mcl{C}$ be a Kan operad which has an h-monoid $(\mu,\eta)$ in the sense of \emph{Definition \ref{def: mon}}.  An \emph{$E_\infty$ structure} on $(\mcl{C},\mu)$ is a map of operads $\mcl{T} \to \mcl{C}$ which maps the generating 2-tree $\sigma$ in $\mcl{T}_2$ to $\mu$.  \end{defi}

There is a space $\mcl{E}_\infty(\mcl{C},\mu)$ of $E_\infty$ structures on $(\mcl{C},\mu)$.  This is a Kan complex in which an $k$-simplex is a map of operads 
$$(\Delta[k] \times \mcl{T}, \Delta[k] \times \sigma) \to (\mcl{C},\mu)\;.$$  The homotopy unit $\eta$ plays no formal r\^{o}le, but its existence is needed so that Lemma \ref{lemma: gastr} can be applied.

\subsection{Stages for an $E_\infty$ structure\label{sec: stage}}

We build up an $E_\infty$ structure stage by stage.  First we define some filtrations on the operad~$\mcl{T}$.

There are two straightforward ways of filtering $\mcl{T}$.  (They are filtrations by subspaces, not by suboperads.)  The first is the \emph{height filtration}: at level $n$, this consists of the spaces $\mcl{T}_S$ with $|S| \le n$.  The other filtration is the \emph{bar filtration}, which is inherited from the filtration of the nerve of $\Iso{\Om}$.   The $n$th stage here is 
$\mcl{B}^n\mcl{T} \; = \; \mcl{P} \times \sk_n\!\mcl{Q}$, where $\sk_n\!\mcl{Q}$ is the $n$-skeleton of the Barratt-Eccles operad.

Attempting to construct $E_\infty$ structures by induction on the height filtration leads to difficulties, as is explained in the introduction of \cite{goe-hop}.  We use instead a filtration which combines the two defined above.

\begin{defi} \label{diag} The \emph{diagonal filtration} on $\mcl{T}$ is defined as the sum of the height filtration and the bar filtration.  Thus, for $n\ge 0$, the $n$th level in the diagonal filtration of $\mcl{T}$ consists of the spaces $\mcl{B}^p\mcl{T}_S$ for $p+|S| \le n$. \end{defi}

\begin{defi} \label{def: stage} An \emph{$n$-stage} $\phi$ for an $E_\infty$ structure on $(\mcl{C},\mu)$ consists of maps 
$\phi_S\colon \mcl{B}^p\mcl{T}_S \to \mcl{C}_S$ for all $p$ and $S$ such that $p + |S| =n$, subject to the following conditions.
\begin{enumerate} \item \emph{Equivariance:} $\phi_S$ commutes with maps induced by isomorphisms of $S$.
\item \emph{Coherence:} the following diagram commutes for every partition $S = M \sqcup N$ (with $|M|,|N^+| \ge 2$)
$$\begin{array}{ccc}
 \mcl{B}^p\mcl{T}_M \times  \mcl{B}^p\mcl{T}_{N^+}  & \stk{\circ_+}{\lra} &  \mcl{B}^p\mcl{T}_S     \\
\phi_M\! \times\! \phi_{N^+} \Big\downarrow \phantom{\phi_M \times \phi_N}     &  & 
\phi_S \Big\downarrow \phantom{\phi_S}                                                                      \\
\mcl{C}_M \times \mcl{C}_{N^+}       &  \stk{\circ_+}{\lra} &  \mcl{C}_S    \end{array} $$ 
where $\circ_+$ is operad composition.
\item  If $\sigma_2$ is the star tree in $\mcl{T}_2$, then $\phi_2(\sigma_2) = \mu$.  \end{enumerate} 
\end{defi}

An $n$-stage defines an $(n-1)$-stage by restriction. A 1-stage specifies a commutativity homotopy and an associativity homotopy for the monoid structure $\mu$. 
The horizontal maps in the above diagram are the restrictions of the composition maps in the operads $\mcl{T}$ and $\mcl{C}$ respectively.   The coherence condition asserts that the $\phi_S$ form a map of operads insofar as this condition is defined.  Indeed, the values of $\phi_S$ on any face of $\mcl{B}^p\mcl{T}_S$ are determined by the $(n-1)$-stage underlying $\phi$, because $|M|, |N| < |S|$.  This is the basis for the inductive procedure for analysing $\mcl{E}_\infty(\mcl{C},\mu)$.

\subsection{A tower for $\mcl{E}_\infty(\mcl{C},\mu)$ \label{sec: tower}}

We denote by $\St^n(\mcl{C},\mu)$ the space of $n$-stages for $E_\infty$ structures on $(\mcl{C},\mu,\eta)$.  (This space is defined by analogy with the space of $E_\infty$ structures described after Definition~\ref{def: formal E} above: a $k$-simplex of this space consists of maps  $\phi_S\colon \Delta[k] \times \mcl{B}^p\mcl{T}_S \to \mcl{C}_S$ satisfying a $\Delta[k]$-parametrized version of the conditions in Definition \ref{def: stage}.)  

\begin{prop}\label{prop: kantower} If $\mcl{C}$ is a Kan operad, then each restriction map $$\St^n(\mcl{C},\mu) \lra \St^{n-1}(\mcl{C},\mu)$$ is a Kan fibration, and its fibre $\St^n_{n-1}$ has the homotopy type of a finite product of  copies of $(n-2)$-fold loop spaces $\Om^{n-2}\mcl{C}_k$ for various~$k$.   Thus there is a tower of Kan fibrations
$$ \dots \lra \St^n(\mcl{C},\mu) \lra \St^{n-1}(\mcl{C},\mu) \lra \dots \lra \St^0(\mcl{C},\mu) \approx \{*\}$$
with limit $\mcl{E}_\infty(\mcl{C},\mu)$. \end{prop}

{\it Proof.} Let us consider the map 
$\St^n(\mcl{C},\mu) \lra \St^{n-1}(\mcl{C},\mu)$.  To extend an $(n-1)$-stage to an $n$-stage, we must extend given maps $\mcl{B}^{p-1}\mcl{T}_S \to \mcl{C}_S$ over $\mcl{B}^p\mcl{T}_S$, naturally with respect to isomorphisms of $S$, for all $p$ and $S$ such that 
$p + |S| = n$.  The extension is in each case already given over $\mcl{B}^p\del\mcl{T}_S$ by the coherence condition in Definition~\ref{def: stage}.  Reducing $\Om$ to a spine,  it suffices to find a $\Sigma_k$-equivariant extension when $S = \underline{k} = \{1,2,\dots,k\}$ for $k = 2,\dots,n$.  This implies that the fibre $\St^n_{n-1}$ of the restriction map 
$\St^n(\mcl{C},\mu) \lra \St^{n-1}(\mcl{C},\mu)$ is
$$\St^n_{n-1} \quad = \quad \prod_{k=2}^n \Map_{\Sigma_k}
\bigl(\mcl{B}^{n-k}\mcl{T}_k/(\mcl{B}^{n-k}\del\mcl{T}_k \cup \mcl{B}^{n-k-1}\mcl{T}_k), \;\mcl{C}_k\bigr)\;.$$  
(Since the inclusion of $\mcl{B}^{n-k}\del\mcl{T}_k \cup \mcl{B}^{n-k-1}\mcl{T}_k$ into $\mcl{B}^{n-k}\mcl{T}_k$ is a cofibration, and $\mcl{C}_k$ is a Kan complex, it follows that this restriction map is a fibration: the equivariance condition is not a problem, because the symmetric group $\Sigma_k$ acts freely on the domain.)

On unravelling the definitions of $\mcl{T}_k$ and $\del\mcl{T}_k$, our fibre becomes 
$$\prod_{k=2}^n \Map_{\Sigma_k}\bigl((\mcl{P}_k/\del\mcl{P}_k) \wedge 
(\sk_{n-k}\!\mcl{Q}_k/\sk_{n-k-1}\!\mcl{Q}_k), \;\mcl{C}_k\bigr)\;. $$
By Proposition~\ref{prop: partbdry} the space $(\mcl{P}_k/\del\mcl{P}_k)$ has the homotopy type of a wedge of \mbox{$(k-2)$}-spheres: its non-trivial homology group is $\mcl{S}_k^*$.  Further, in $\mcl{Q}_k$ the quotient of skeleta  
$\sk_{n-k}\!\mcl{Q}_k/\sk_{n-k-1}\!\mcl{Q}_k$ is a wedge of $(n-k)$-spheres, indexed by $(n-k)$-tuples of elements of $\Sigma_k$.  The action of $\Sigma_k$  on this wedge is free.  It follows that the fibre $\St^n_{n-1}$ is a product of copies of $\Om^{n-2}\mcl{C}_k$ for various~$k$ with $2\le k \le n$.  

Since the base of the tower 
$$ \dots \lra \St^n(\mcl{C},\mu) \lra \St^{n-1}(\mcl{C},\mu) \lra \dots \lra \St^0(\mcl{C},\mu) \approx \{*\}$$ is a point, and all the maps are fibrations, it follows by recursion that all spaces in the tower are Kan complexes.  The limit of the tower is evidently the space $\mcl{E}_\infty(\mcl{C},\mu)$ of $E_\infty$ structures. $\Box$

\subsection{The homotopy spectral sequence\label{sec: hoss}}

The homotopy exact sequences of the fibrations 
$$\St_{p-1}^p \lra \St^p(\mcl{C},\mu) \lra \St^{p-1}(\mcl{C},\mu)$$
in the tower of \ref{sec: tower} yield an unwound exact couple 
and therefore a homotopy spectral sequence for $\pi_*\mcl{E}_\infty(\mcl{C},\mu)$ in the form 
$$E_1^{p,q} \approx \pi_{q-p}\St_{p-1}^p \Longrightarrow \pi_{q-p}\,\mcl{E}_\infty(\mcl{C},\mu)\;.$$  

The properties of the homotopy spectral sequence of a tower of fibrations are treated comprehensively in (\cite{goe-jar}, Ch.\ VI), which is our principal reference on this topic.  Two particular problems arise.

The minor problem is that the exact sequences above end in maps $\pi_0\St^p(\mcl{C},\mu) \lra \pi_0\St^{p-1}(\mcl{C},\mu)$ which need not be surjective.  This means that the spectral sequence has a \emph{fringe} at total degree 0 rather than an edge, and a separate investigation is needed to determine which elements on the fringe actually arise from elements of $\pi_0\mcl{E}_\infty(\mcl{C},\mu)$.  We shall associate an obstruction theory with this fringe, like the one developed in detail for an analogous situation by \cite{bousfield:1989}.  

The major problem with the spectral sequence is that convergence is not automatic, and requires separate arguments.   We discuss this further in \S\ref{sec: conv}. 

\subsection{The $E_1$ and $E_2$ terms of the spectral sequence\label{sec: e1}} 

We shall usually assume that  $(\mu,\eta)$ is a loop h-monoid in $\mcl{C}$ in the sense of Definition \ref{looperad}.

We know from \S\ref{sec: tower} and \S\ref{sec: hoss} that
$$\begin{aligned} E_1^{p,q} \quad \approx& \quad \pi_{q-p}\St_{p-1}^p \\
            \approx& \quad  \pi_{q-p}\biggl( \prod_{k=2}^p \Map_{\Sigma_k}
\bigl(\mcl{B}^{p-k}\mcl{T}_k/(\mcl{B}^{p-k}\del\mcl{T}_k \cup \mcl{B}^{p-k-1}\mcl{T}_k), \;\mcl{C}_k\bigr)\biggr)\; . \end{aligned}$$
Since $\mcl{B}^{p-k}\mcl{T}_k/(\mcl{B}^{p-k}\del\mcl{T}_k \,\cup \,\mcl{B}^{p-k-1}\mcl{T}_k)$ is equivariantly homotopy equivalent to some wedge of $(p-2)$-spheres with free $\Sigma_k$-action, this implies that 
$$\begin{aligned} E_1^{p,q} \;\; \approx& \;\; \bigoplus_{k=2}^p \;
\Hom_{\Sigma_k}\bigl(\til H_{p-2}(\mcl{B}^{p-k}\mcl{T}_k/(\mcl{B}^{p-k}\del\mcl{T}_k \,\cup \, \mcl{B}^{p-k-1}\mcl{T}_k)), \;\pi_{q-2}\,\mcl{C}_k\bigr)   \\
\approx& \;\; \bigoplus_{k=2}^p \;
\til H^{p-2}_{\Sigma_k}\bigl(\mcl{B}^{p-k}\mcl{T}_k/(\mcl{B}^{p-k}\del\mcl{T}_k \,\cup \, \mcl{B}^{p-k-1}\mcl{T}_k); \;\pi_{q-2}\,\mcl{C}_k\bigr) \, .\end{aligned}$$
This description of $E_1^{p,q}$ in terms of equivariant cohomology arises from Proposition \ref{prop: partbdry}.  In order to calculate $d_1^{p,q}$ however we need to look at a chain-level description: we must write elements of $E_1^{p,q}$ as relative cocycles of pairs
$(\mcl{B}^{p-k}\mcl{T}_k, \mcl{B}^{p-k}\del\mcl{T}_k \cup \mcl{B}^{p-k-1}\mcl{T}_k)$ and calculate the connecting coboundaries.

The operad space $\mcl{T}_S = \mcl{P}_S \times \mcl{Q}_S$ has a cellular structure in which the cells are products $\kappa \times \psi$, where $\kappa$ is a simplex of the partition operad space $\mcl{P}_S$ and $\psi$ a simplex of $\mcl{Q}_S$.
By using the chain complex $C_*(\mcl{T}_S)$ of this prismatic cell structure, we avoid the Eilenberg-Zilber theorem: we  have an isomorphism 
$$C_*(\mcl{T}_S) \quad \approx \quad C_*(\mcl{P}_S) \otimes C_*(\mcl{Q}_S)\;.$$
We compare $C_*(\mcl{T}_S)$ with the chains on the nerve $\til N(S/\Om)$.
Suppose that $[f_1|f_2|\dots|f_q]$ is a simplex of $\til N(S/\Om)$ having width $p \ge 0$. Then the face 
$$S_{q-p} \stk{f_{q-p+1}}{\lla} \dots \stk{f_{q}}{\lla} S_q \stk{\approx}{\lla} S$$
is a $p$-simplex of the Barratt-Eccles space $\mcl{Q}_S$, and the face 
$$S_0 \stk{f_1}{\lla}S_1 \stk{f_2}{\lla} \dots \stk{f_{q-p-1}}{\lla} S_{q-p-1} \stk{\eta}{\lla} S$$ (where $\eta$ is the composite map) is
a non-degenerate $(q-p-1)$-simplex of the nerve $N(S/\!/\Om)$, which by Lemma \ref{lemma: slicetree} of \S\ref{sec: partree} is related to $\mcl{P}_S$ by a homotopy equivalence of pairs under which $N(S/\!/\Om/\!/\ulin{1})$ corresponds to $\del\mcl{P}_S$.  The chain map 
$$[f_1|f_2|f_3|\dots|f_q] \quad \mapsto \quad [f_1|f_2|\dots|f_{q-p-1}] \otimes 
[f_{q-p+1}|\dots|f_q]$$ is now a chain isomorphism of degree~-1
$$\begin{aligned} C_*(\til N(S/\Om), \til N(S/\Om/\!/\ulin{1})) \quad &\approx \quad C_*(\mcl{P}_S, \del\mcl{P}_S) \otimes C_*(\mcl{Q}_S) \\
\quad &\approx \quad C_*(\mcl{T}_S, \del\mcl{T}_S)\;. \end{aligned}$$ 
(The map $f_{q-p}$, having become structural in $\til N\!(S/\!/\Om)$, no longer appears.)
The resulting equivalence of pairs is an equivalence of filtered pairs, since the width filtration on the left corresponds precisely to the bar filtration (on $\mcl{Q}_S$ and $\mcl{T}_S$) on the right. We may therefore pass to filtration quotients.  Denoting as before by $\mcl{W}^p$ the $p$th width-filtration stage, and by $\mcl{W}^p/\mcl{W}^{p-1}$ the $p$th quotient, we obtain 
$S$-natural isomorphisms 
$$\begin{aligned}
\til C_*\bigl((\mcl{W}^p/\mcl{W}^{p-1})(\til N(S/\Om,\,S/\Om/\!/\ulin{1})\bigr) \quad \approx& 
\quad \til C_*\bigl((\mcl{B}^p/\mcl{B}^{p-1})(\mcl{T}_S,\,\del\mcl{T}_S)\bigr) \\
=& \quad \til
C_*\bigl(\mcl{B}^p\mcl{T}_S/(\mcl{B}^p\del\mcl{T}_S \cup \mcl{B}^{p-1}\mcl{T}_S)\bigr) \,. \end{aligned} $$
When we restrict to the spine of $\Om$ and set $S = \ulin{k}$, the $S$-naturality becomes $\Sigma_k$-equivariance, and the isomorphisms in $\til N(S/\Om)$ become the elements of $\Sigma_k$.  Therefore we can write equivariant cochains on $\mcl{B}^p\mcl{T}_k$ as a bar construction module over $\Sigma_k$
 $$\begin{aligned} \til C_{\Sigma_k}^*\bigl(\mcl{B}^p\mcl{T}_k/(\mcl{B}^p\del\mcl{T}_k \cup \mcl{B}^{p-1}&\mcl{T}_k);
\,\pi_*\,\mcl{C}_k\bigr)   \\ 
 \quad \approx& \quad \mcl{B}^p\bigl(C_*\bigl(N(\ulin{k}/\!/\Om),\, N(\ulin{k}/\!/\Om/\!/\ulin{1})\bigr), \,\Sigma_k,\,\pi_*\,\mcl{C}_k\bigr)\end{aligned}$$ where on the right $\mcl{B}^p(-, \Sigma_k, \pi_*\mcl{C}_k)$ denotes $\Hom(- \otimes (\zed\Sigma_k)^{\otimes p}, \, \pi_*\mcl{C}_k)$.  The differential in the relative cochains on the left corresponds to that in $C_*(N(\ulin{k}/\!/\Om))$ on the right, so by Lemma~\ref{lemma: slicetree} and Proposition~\ref{prop: partbdry}
$$ \til H_{\Sigma_k}^{p+k-2}\bigl(\mcl{B}^p\mcl{T}_k/(\mcl{B}^p\del\mcl{T}_k \cup \mcl{B}^{p-1}\mcl{T}_k);
\,\pi_*\,\mcl{C}_k\bigr)   \quad \approx \quad \mcl{B}^p(\mcl{S}_k^*, \Sigma_k, \pi_*\mcl{C}_k)$$ where as before $\mcl{S}_k^*$ is the superlie representation of Definition~\ref{def: super}; and the groups $\til H_{\Sigma_k}^{j}$ are zero for all $j \ne p+k-2$.  Using our formula for $E_1^{p,q}$ in terms of equivariant homology, we therefore have 
$$\begin{aligned} E_1^{p,q} \;\;\approx& \;\; \bigoplus_{k=2}^p \;
\til H^{p-2}_{\Sigma_k}\bigl(\mcl{B}^{p-k}\mcl{T}_k/(\mcl{B}^{p-k}\del\mcl{T}_k \,\cup \, \mcl{B}^{p-k-1}\mcl{T}_k); \;\pi_{q-2}\,\mcl{C}_k\bigr) \\
\approx& \;\; \bigoplus_{k=2}^p \; \mcl{B}^{p-k}(\mcl{S}_k^*, \Sigma_k, \pi_{q-2}\mcl{C}_k)
\, .\end{aligned}$$
Now we need to calculate $d_1^{p,q}$ in these terms.

\begin{lemma}\label{lemma: d1} Commutativity holds in the diagram 
$$\begin{array}{ccccc}
E_1^{p,q}  &  &{\approx}  & &
\bigoplus_{k=2}^p \; \mcl{B}^{p-k}(\mcl{S}_k^*, \Sigma_k, \pi_{q-2}\mcl{C}_k)  \\  &&&& \\
d_1^{p,q}\Big\downarrow  &&&& \Big\downarrow d_1^{p,-1}  \\  &&&& \\
E_1^{p+1,q}  & & \approx  & &
\bigoplus_{k=2}^{p+1} \; \mcl{B}^{p+1-k}(\mcl{S}_k^*, \Sigma_k, \pi_{q-2}\mcl{C}_k) 
\end{array} $$
in which the horizontal isomorphisms are those constructed above.  On~the left, 
$d_1^{p,q}$ is the differential in the homotopy spectral sequence; on~the right,  $d_1^{p,-1}$ is the algebraically-defined morphism in  Corollary~\ref{cor: finalmess}.
\end{lemma} 

{\it Proof.}  On the left,  we have the geometry of operads: the differential $d_1^{p,q}$ arises from the face structure in the tree operad $\mcl{T}$.  On the right, we have the algebra of the diagonally-filtered complex $\Rupga^*(\pi_*\mcl{C})$: the differential $d_1^{p,-1}$  is defined in terms of the face operators in the nerve $\til N(S/\Om)$.
The connection between the modules on the left and right arises through the isomorphism of degree~1 proved above  
$$C_*(\mcl{T}_S) \quad =  \quad C_*(\mcl{P}_S) \otimes C_*(\mcl{Q}_S) \quad 
\approx \quad C_*(\til N(S/\Om))\;.$$ 
Consider a cochain $\gamma$ in a single summand 
$\mcl{B}^{p-k}(\mcl{S}_k^*, \Sigma_k, \pi_{q-2}\mcl{C}_k)$ in the top right.  The image (or coboundary) $d_1^{p,-1}(\gamma)$ in the lower right can take non-trivial values only on linear combinations of simplices
$$S_1 \stk{f_2}{\lla} \dots \stk{f_{q-p}}{\lla} \dots \stk{f_{q}}{\lla} S_q \stk{\approx}{\lla} S$$
having width $p-k$ or $p-k+1$.  The coboundary $d_1^{p,-1}(\gamma)$ is the alternating sum of coface operators $\delta^i$ dual to the face operators in $\Rdnga_*(\pi_*\mcl{C})$.  The sum $\sum_{i=1}^{k-1}(-1)^i\delta^i(\gamma)$ is zero, because $\gamma$ has survived to the $E_1$ term of the algebraic spectral sequence.  There remain $\delta^0(\gamma)$, which is dual to the face $\del_0$ in $\Rdnga_*(\pi_*\mcl{C})$; the cofaces $\delta^i(\gamma)$ for $k \le i \le q$, which are the usual cofaces in the bar construction over $\Sigma_k$; and $\delta^{q+1}(\gamma)$, which is the external final coface.

Now we look at the element $\olin{\gamma} \in E^1_{p,q}$ which corresponds to $\gamma$ under the horizontal isomorphism.  From \S\ref{sec: e1} we know that $\olin{\gamma}$ is represented by a $\Sigma_k$-equivariant map 
$$\mcl{B}^{p-k}\mcl{T}_k \; \lra \; \Om^{q-p}\mcl{C}_k$$
which is trivial on the boundary $\mcl{B}^{p-k}\del\mcl{T}_k \cup \mcl{B}^{p-k-1}\mcl{T}_k$. The coboundary $d_1^{p,q}(\olin{\gamma})$ has components in two summands of $E^1_{p+1,q}$, namely in $\mcl{B}^{p+1-k}(\mcl{S}_k^*, \Sigma_k, \pi_{q-2}\mcl{C}_k)$ and in $\mcl{B}^{p-k}(\mcl{S}_{k+1}^*, \Sigma_{k+1}, \pi_{q-2}\mcl{C}_{k+1})$.   The required elements in each are determined by the inclusions of $\mcl{B}^{p-k}\mcl{T}_k$ as a face in respectively either $\mcl{B}^{p+1-k}\mcl{T}_k$ or $\mcl{B}^{p-k}\mcl{T}_{k+1}$.  In the first case we have the usual face inclusions in the bar construction, which correspond to the cofaces $\delta^i(\gamma)$ for $k \le i \le q$ in the algebraic situation.  In the second case the face inclusions are induced by the operad structure in $\mcl{P}$: the face inclusions 
$$\mcl{P}_k \;\approx\;\mcl{P}_k \times \mcl{P}_2 \; \lra \; \mcl{P}_{k+1}$$ (which are all conjugate under the symmetric group) correspond to the algebraic coface $\delta^0$, and the inclusions $$\mcl{P}_k \;\approx\;\mcl{P}_2 \times \mcl{P}_k \; \lra \; \mcl{P}_{k+1}$$ correspond to the remaining algebraic coface $\delta^{q+1}$.  (The remaining faces, the inclusions of $\mcl{P}_u \times \mcl{P}_v$ with $u,v > 2$, contribute zero in the $E^1$ term: they give rise to higher differentials.)

The identification  of the algebraic and geometric differentials establishes the Lemma. $\Box$
\medskip

Now for the main result.  We assume the data of Lemma~\ref{lemma: gastr}: $\mcl{C}$ is a Kan operad which has a loop h-monoid $(\mu,\eta)$ where $\mu\in \mcl{C}_2$ and $\eta \in \mcl{C}_0$, and we seek information about homotopy groups of the space $\mcl{E}_\infty(\mcl{C},\mu)$ of $E_\infty$ structures on $\mcl{C}$ which extend  $\mu$.  (This space may be empty.)

\begin{thm}\label{thm: e2} The homotopy spectral sequence for $E_\infty$ structures on an operad $\mcl{C}$ with a loop h-monoid is fringed at total degree $p-q = 0$: it has 
$$E_2^{p,q} \;\; \approx \;\; \pi^{p-1}\pi_{q-2}\,\mcl{C} \quad \Longrightarrow \quad 
\pi_{q-p}\,\mcl{E}_\infty(\mcl{C},\mu)$$
$$d_r^{p,q}\colon E_r^{p,q} \lra E_r^{p+r, q+r-1}$$
where $\pi_*\,\mcl{C}$ is the homotopy $\Ga$-module of Lemma~\ref{lemma: gastr} and $\pi^*$ denotes the stable cohomotopy of a $\Ga$-module. \end{thm}

{\it Proof.} This is a standard Bousfield-Kan fringed homotopy spectral sequence for a tower of fibrations (\cite{bous-kan};~\cite{goe-jar}, Ch.\,VI \S2).  By Lemma~\ref{lemma: d1} above, the  cohomology of $(E_1, d_1)$ is identified with the cohomology of the complex in  Corollary~\ref{cor: finalmess} of \ref{sec: redcx}.  This gives as claimed  
$$E_2^{p,q} \;\; \approx \;\; \pi^{p-1}\pi_{q-2}\,{\mcl{C}}\;.$$  The abutment is as usual the homotopy $\pi_{*}\,\mcl{E}_\infty(\mcl{C},\mu)$ of the limit of the tower of fibrations.  Conditions for the somewhat delicate convergence of the spectral sequence are investigated in~\ref{sec: conv} and \ref{sec: fringe} below.  $\Box$
\smallskip

\subsection{Conditions for convergence\label{sec: conv}}

Since the tower of fibrations is bounded below with $\St^0(\mcl{C},\mu) \approx \{*\}$, each term $E_1^{p,q}$ is modified by only finitely many incoming differentials $d_1, \dots, d_p$.  Therefore $E_{s+1}^{p,q} \subset E_s^{p.q}$ for all $s>p$.  We set 
$$E_\infty^{p,q} \quad = \quad \lim_s  E_s^{p,q} \quad \approx \quad \bigcap_{s>p}  E_s^{p,q}$$
and we note that the derived limit $\lim_s^1  E_s^{p,q}$ is also defined.  

The result we should like to have is \emph{complete convergence} in the terminology of \cite{goe-jar}, which is equivalent to the validity of the two following statements:

\begin{enumerate}
\item For all $q-p > 0$  (that is, away from the fringe) the $E_\infty^{p,q}$ are the quotients in a descending composition series for the group 
$$\lim_s \pi_{q-p}\St^s(\mcl{C},\mu)\;. $$

\item The natural map 
$$ \pi_*\mcl{E}_\infty(\mcl{C},\mu)  \; = \;  \pi_*(\lim_s \St^s(\mcl{C},\mu))
                                                                   \lra \lim_s \pi_*\St^s(\mcl{C},\mu) $$
is an isomorphism in all degrees.
\end{enumerate}

Either of these conditions may fail for the homotopy spectral sequence of a general tower of fibrations.  Some composition factors may for instance be proper subgroups of the corresponding $E_\infty^{p,q}$. 

There is a standard condition for complete convergence.  We refer to 
(\cite{goe-jar} Ch.\ VI Lemma 2.20) for the proof of the following.

\begin{prop}\label{prop: lim1} The spectral sequence of {Theorem \ref{thm: e2}} is completely convergent if and only if the derived limit $\lim_s^1 E_s^{p,q}$ is trivial for all $q-p > 0$. $\Box$ \end{prop}

One may hope to satisfy the hypothesis here by showing that the systems $\{E_s^{p,q}\}_s$ satisfy a Mittag-Leffler condition. 

\begin{cor}\label{cor: fin} If all the homotopy groups $\pi_j\mcl{C}_k$ of the operad spaces are finite, then the spectral sequence is completely convergent. \end{cor}

{\it Proof.}  In this case Lemma~\ref{lemma: d1} shows that all the $E_1^{p,q}$ groups are finite, so that the Mittag-Leffler condition guarantees the vanishing of the derived limits required by Proposition~\ref{prop: lim1}. $\Box$

It is possible to generalize this Corollary to a situation where the spectral sequence were in the category of modules over a ring, and the $\pi_j\mcl{C}_k$ were all modules of finite length.

\subsection{Obstructions and the fringe\label{sec: fringe}}

The spectral sequence can be used as normal to calculate $\pi_i \mcl{E}_\infty(\mcl{C},\mu)$ for $i \ge 1$.  The calculation of $\pi_0 \mcl{E}_\infty(\mcl{C},\mu)$ (in particular, the question of whether any $E_\infty$ structure exists) encounters the fringe.  Explicitly, the spectral sequence arises from the homotopy exact sequences of the fibrations $\St^p(\mcl{C},\mu) \to \St^{p-1}(\mcl{C},\mu)$; and such a sequence ends 
$$E_1^{p,p} = \pi_0\St_{p-1}^p \; \lra \; \pi_0\St^p(\mcl{C},\mu) \; \lra \; \pi_0\St^{p-1}(\mcl{C},\mu)$$
with no indication of whether the last map is surjective.  This fringe in total degree zero is a normal feature of the homotopy spectral sequence of a tower.

The tower for $\mcl{E}_\infty(\mcl{C},\mu)$ has extra structure.  Since the tree space $\mcl{P}_k$ is the cone on $\del\mcl{P}_k$ by Proposition~\ref{prop: partop}, it follows that the bar filtration stage $\mcl{B}^n\mcl{T}_k$ (defined in~\ref{sec: stage}) of the tree operad is obtained from $\mcl{B}^{n-1}\mcl{T}_k$ by attaching $k^n$ copies of 
$\mcl{P}_k \times \Delta^n$ along their boundaries
$\del\mcl{P}_k \times \Delta^n \, \cup \, \mcl{P}_k \times \del\Delta^n$.  These are  $\Sigma_k$-equivariant \emph{principal} cofibrations.  Using the attaching maps we can therefore extend each cofibration sequence  $$\mcl{B}^{n-1}\mcl{T}_k  \; \subset \; \mcl{B}^{n}\mcl{T}_k$$ one step to the left.

Passing to equivariant mapping spaces of the diagonal filtration, we have shown that 
each $\St^p(\mcl{C},\mu) \to \St^{p-1}(\mcl{C},\mu)$ is a principal fibration.  We have de-looped the fibre $\St_{p-1}^p$ and we can extend the homotopy exact couple by one step to the right, so that its ends
$$\pi_0\St^p(\mcl{C},\mu) \; \lra \; \pi_0\St^{p-1}(\mcl{C},\mu) \; \lra \; E_1^{p,p-1} \approx \bigoplus_{k=2}^p \; \mcl{B}^{p-k}(\mcl{S}_k^*, \Sigma_k, \pi_{p-3}\mcl{C}_k)$$ where we have identified the new term $E_1^{p,p-1}$ as in \S\ref{sec: e1}.  Therefore we have in the group on the right an obstruction in the cochain complex of (Corollary~\ref{cor: finalmess}, \S\ref{sec: redcx}) to lifting a $(p-1)$-stage for a structure to a $p$-stage.  As usual, this cochain must be a cocycle, as its value on a chain is calculated on a relative boundary.

The earlier part of the homotopy exact sequence implies that there are difference chains lying in $E_1^{p-1,p-1}$ for $(p-1)$-stages extending the same underlying $(p-2)$-stage; and all the chains in this group occur as differences.  As in Eilenberg's original theory we may ask the effect on the obstruction cocycle of altering the $(p-1)$-stage by a general difference cochain. This amounts to asking for an evaluation of the image of the differential 
$d_1^{p-1,p-1}\colon E_1^{p-1,p-1} \to E_1^{p,p-1}$.  The answer is given by the desuspension of Lemma~\ref{lemma: d1}: one can alter the cocycle by any coboundary in the algebraic complex of Corollary~\ref{cor: finalmess}.  Therefore we have a primary obstruction in $E_2^{p,p-1} \approx \pi^{p-1}\pi_{p-3}\mcl{C}$ to replacing the $(p-1)$-stage by a $p$-stage with the same underlying $(p-2)$-stage.  (We have no description of the differential $d_2$ on the extended fringe $E_2^{p,p-1}$, and therefore do not define $E_r^{p,p-1}$ for $r>2$.) This discussion proves the following theorem.

\begin{thm}~\label{thm: obst} Let $\mcl{C}$ be an operad of Kan spaces, and let $\mu \in \mcl{C}_2$ be a loop h-monoid in $\mcl{C}$.
\begin{enumerate} \item 
Given, for some $n \ge 3$, an $n$-stage of an $E_\infty$ structure which extends $\mu$, there is an obstruction in $ \pi^{n}\pi_{n-2}\,\mcl{C}$ which vanishes if and only if the underlying $(n-1)$-stage extends to a $(n+1)$-stage.
\item If $\pi^{n}\pi_{n-2}\,\mcl{C} \approx 0$ for all $n \ge 3$, then there exists an $E_\infty$ structure extending the homotopy-commutative, homotopy-associative multiplication $\mu$.  
\item If also 
$\pi^{n}\pi_{n-1}\mcl{C} \approx 0$ for all $n \ge 2$, then every finite stage of this $E_\infty$ structure is unique up to homotopy.   If further the spectral sequence is completely convergent in the sense of~\ref{sec: conv}, then the $E_\infty$ structure extending $\mu$ is unique up to homotopy.~$\Box$ \end{enumerate} \end{thm}

\subsection{Application: $E_\infty$ structures on ring spectra\label{sec: spectra}}

As a special case of Theorem~\ref{thm: obst} we derive a generalization of Theorems 5.5 and 5.6 of~\cite{rob:inv}.

We work with spectra in a simplicial model category $\mcl{S}$ in which the smash product is symmetric monoidal, and the homotopy category $\Ho(\mcl{S})$ is equivalent to the standard stable homotopy category.  Thus either the category of orthogonal spectra \cite{mmss} or the category of symmetric spectra \cite{hovess} for $\mcl{S}$ would be suitable.

Let $V$ be a ring spectrum in the classical sense: thus we are given a multiplication 
$\mu\colon V \wedge V \to V$ and a unit map $\eta\colon S \to V$, making $V$ into a commutative monoid in the homotopy category.  Then $(\mu,\eta)$ is an h-monoid structure in the endomorphism operad $\Calend(V)$, where 
$$\Calend(V)_n \quad = \quad \Map_{\mcl{S}}(V^{\wedge n}, V)\;;$$
and the invertibility of suspension in $\mcl{S}$ ensures that $(\mu,\eta)$ is a loop h-monoid structure in the sense of~\ref{sec: piop}.  An $E_\infty$ structure on $V$ is simply a morphism of operads $\mcl{T} \to \Calend(V)$; and we may use the theorems of \S\S\ref{sec: modsp}--\ref{sec: fringe} to investigate these.  The coefficients for the obstruction theory are in right $\Ga$-modules
$$\begin{aligned} \pi_{n-i}\,\Calend(V)_* \quad 
&= \quad \pi_{n-i}\Map_{\mcl{S}}(V^{\wedge *}, V) \\
&= \quad V^{i-n}(V^{\wedge*}) \end{aligned}$$ 
where we perhaps should state explicitly that for each value of $n-i$ these groups form a right $\Ga$-module by Lemma~\ref{lemma: gastr}. 
From Theorem~\ref{thm: obst} we deduce

\begin{thm}\label{thm: inf}  Let $V$ be a ring spectrum with homotopy commutative, homotopy associative multiplication 
$\mu\colon V \wedge V \to V$ and unit $\eta\colon S \to V$.  
\begin{enumerate} \item 
Given, for some $n \ge 3$, an $n$-stage of an $E_\infty$ structure on $V$ extending $\mu$, there is an obstruction in $\pi^{n}\bigl(V^{2-n}(V^{\wedge*})\bigr)$ which vanishes if and only if the underlying $(n-1)$-stage extends to a $(n+1)$-stage.
\item If $\pi^{n}\bigl(V^{2-n}(V^{\wedge*})\bigr) \approx 0$ for all $n \ge 3$, then there exists an $E_\infty$ structure extending the homotopy-commutative, homotopy-associative multiplication $\mu$.  
\item If also 
$\pi^{n}\bigl(V^{1-n}(V^{\wedge*})\bigr) \approx 0$ for all $n \ge 2$, then every finite stage of this $E_\infty$ structure is unique up to homotopy.   If further the spectral sequence of~\ref{sec: hoss} is completely convergent in the sense of~\ref{sec: conv}, then the $E_\infty$ structure extending $\mu$ is unique up to homotopy.~$\Box$ \end{enumerate} \end{thm}

We recall that the convergence condition in the last clause of Theorem~\ref{thm: inf} is non-trivial.   By Corollary~\ref{cor: fin} a sufficient condition for complete convergence is that all groups $V^*(V^{\wedge*})$ be finite.  This suggests that if only finite type conditions are available, it would be worth looking at completions.

A frequently-arising special case is when $V$ satisfies a universal coefficient theorem in the form
$$V^*(V^{\wedge k}) \; \approx \; \Hom_{V_*}\bigl((V_*V)^{\otimes k}, V_*\bigr)$$
which is not needed in Theorem~\ref{thm: inf}, but is true whenever the bialgebra $V_*V$ of homology co-operations is projective over the coefficient ring $V_*$.  Then the graded right $\Ga$-module $\pi_*\,\!\Calend(V)$ is the dual of the graded Loday functor $\mcl{L}(V^*V)$ of (\cite{pir:hodge}, \S1.7), so  that $\pi^*\!\pi_*\,\!\Calend(V)$ is the $\Ga$-cohomology~\cite{rob-whi2} of $V^*V$.  Then from Theorem~\ref{thm: inf} one recovers Theorems 5.5 and 5.6 of~\cite{rob:inv}, corrected by the insertion of a convergence clause which was inadvertently omitted from 5.6 of~\cite{rob:inv}.

\bibliography{math}
\bibliographystyle{amsalpha}
\bigskip
\end{document}